\documentclass[10pt,twoside, a4paper]{article}
\usepackage[utf8]{inputenc}
\usepackage[T1]{fontenc}
\usepackage{amsfonts, amsmath}
\usepackage{txfonts}
\usepackage{geometry}
\usepackage[calcwidth,  bf,  compact]{titlesec}
\usepackage[titletoc]{appendix}
\usepackage{etoolbox, apptools}
\AtBeginEnvironment{appendices}{\appendixtrue}
\usepackage{setspace}
\titleformat{\part}[block]{\bf }
{\large}{0pt}{\large}
\titleformat{\subsubsection}[block]{\large\bf}
{\thesubsubsection. }{3pt}{\raggedright}
\titleformat{\subsection}[block]{\bf}
{\raggedright 
\thesubsection. }{1pt}{
\raggedright}
\titleformat{\section}[block]{\bf}
{\large \IfAppendix{\appendixname~}{\relax}\thesection. \IfAppendix{: }{}}{0pt}{\large}
 
\usepackage{amssymb}
\usepackage{mathtools}
\usepackage[mathscr]{eucal}
\usepackage{booktabs}

\def\bibfont{\footnotesize}
\usepackage[natural, dvipsnames]{xcolor}
\usepackage{savesym}
\savesymbol{AND}
\usepackage{algorithm}
\usepackage[noend]{algpseudocode}
\usepackage[amsmath, thref,  thmmarks,  hyperref]{ntheorem}
\usepackage{nameref}
\usepackage[pdftex, hidelinks]{hyperref}
\usepackage[square, numbers]{natbib}
\usepackage{cleveref}

\newtheorem{theorem}{Theorem}

\newtheorem{lemma}{Lemma}

\newtheorem{corollary}{Corollary}

\newtheorem{proposition}{Proposition}

\newtheorem{definition}{Definition}

\newtheorem{remark}{Remark}

\font\QEDlogofont=msam10 at 10pt
\def\QEDlogo{\hbox{\QEDlogofont\char'003}}

\theoremsymbol{\QEDlogo}
\theoremheaderfont{\hskip\normalfont\itshape}
\theorembodyfont{\normalfont}
\theoremseparator{.}
\theoremstyle{nonumberplain}

\definecolor{labelkey}{rgb}{0.6, 0, 1}
\newcommand{\e}{\mathrm E}
\newcommand{\E}[2][]{{\mathrm E}_{#1}\left(#2\right)}

\newcommand{\I}[1]{{\mathrm I}_{#1}}

\newcommand{\R}{\mathbb{R}}
\newcommand{\Si}{\mathbb{S}}

\newcommand{\N}{{\mathbb N}}
\newcommand{\bs}[1]{\boldsymbol {#1}}

\newcommand{\pfrac}[2]{\left(\frac{#1}{#2}\right)}

\renewcommand{\d}{\mathrm{d}}
\newcommand{\Be}{\mathrm{B}}

\newcommand{\indep}{\protect\mathpalette{\protect\independenT}{\perp}}
\DeclareMathOperator*{\argmax}{\mathrm{argmax}}
\def\independenT#1#2{\mathrel{\rlap{$#1#2$}\mkern2mu{#1#2}}}
\newlength{\depthofsumsign}
\newcommand{\eqdis}{\stackrel{\mbox{\tiny d}}{=}}
\newcommand{\comp}{{\mathsf{c}}}   

\crefformat{equation}{Eq.~(#2#1#3)}

\def\d{{\,\textrm{d}}}
\def\and{{\small\textsc{and }\;}}
\newcommand{\email}[1]{\href{mailto:#1}{\texttt{#1}}}

\makeatletter
\def\maketitle{%
	\null
	\thispagestyle{empty}%
	\par
	\begin{center}
		  \Large \textbf{\@title }
	\end{center}
	\par
	\begin{center}
		\normalfont \@author
	\end{center}
	\null
}
\makeatother
\title{Extremal attractors of Liouville copulas}

\author{L\'eo R. Belzile\footnote{\'Ecole Polytechnique F\'ed\'erale de 
Lausanne, EPFL-SB-MATH-STAT,
Station 8, CH-1015
Lausanne, Switzerland. \email{leo.belzile@epfl.ch}} \and Johanna G. 
Ne\v{s}lehov\'{a}\footnote{Department of Mathematics and 
Statistics, 
McGill University, 
805 rue Sherbrooke Ouest, 
Montr\'eal, Qu\'ebec, 
H3A 0B9, Canada. \email{johanna.neslehova@mcgill.ca}}}
\date{}

\usepackage{fancyhdr}
\begin{document}
\hypersetup{pageanchor=false}
\maketitle

\noindent 
This is a copyedited, author-produced version of an article accepted for publication following peer-review in the \textsl{Journal 
of Multivariate Analysis}, an Elsevier publication, \copyright$\,$ 2017. This manuscript version is made available under the 
\href{http://creativecommons.org/licenses/by-nc-nd/4.0/}{CC-BY-NC-ND 4.0 license}. Please cite as 

\begin{quote}
Belzile, L. R and J. G. Ne\v{s}lehov\'{a}. Extremal attractors of Liouville copulas, \textsl{Journal of Multivariate 
Analysis} (2017), 160C, pp. 68--92. \href{https://doi.org/10.1016/j.jmva.2017.05.008}{\texttt{doi:10.1016/j.jmva.2017.05.008}}
\end{quote}

\begin{abstract}
Liouville copulas introduced in \cite{McNeil/Neslehova:2010} are asymmetric 
generalizations  of the ubiquitous Archimedean copula class. They are the dependence 
structures of scale mixtures of Dirichlet distributions, also called Liouville 
distributions. In this paper, the limiting extreme-value attractors of Liouville copulas and 
of their survival counterparts are derived. The limiting max-stable models, termed here 
the scaled extremal Dirichlet, are new and encompass several existing classes of 
multivariate max-stable distributions, including the logistic, negative logistic and 
extremal Dirichlet. As shown herein, the stable tail dependence function and angular 
density of the scaled extremal Dirichlet model have a tractable form, which in turn leads 
to a simple de Haan representation. The latter is used to design efficient algorithms for 
unconditional simulation based on the work of \cite{Dombry:2016} and to derive tractable 
formulas for maximum-likelihood inference. The scaled extremal 
Dirichlet model is illustrated on river flow data of the river Isar in southern Germany. 
\end{abstract}
%

\hypersetup{pageanchor=true}

\section{Introduction}

Copula models play an important role in  the analysis of multivariate data and find 
applications in many areas, including biostatistics, environmental sciences, finance, 
insurance, and risk management. The popularity of copulas is rooted in the 
decomposition of Sklar \cite{Sklar:1959}, which is at the heart of flexible statistical models and various measures, concepts 
and   orderings of dependence between random variables. According to Sklar's result, the distribution function of any random 
vector   $\bs{X} =(X_1, \ldots, X_d)$ with continuous univariate margins $F_1, \dots, F_d$ satisfies, for any $x_1, 
\dots, x_d \in 
\mathbb{R}$, 
\begin{equation*}
\Pr (X_1 \le x_1, \ldots, X_d \le x_d) = C \{ F_1 (x_1), \ldots, F_d (x_d) \}, 
\end{equation*}
for a unique copula $C$, i.e., a distribution function on 
$[0, 1]^d$ whose univariate margins are standard uniform. Alternatively, 
Sklar's decomposition also holds for survival functions, i.e., for any $x_1, \dots, x_d \in \mathbb{R}$, 
\begin{equation*}
\Pr (X_1 > x_1, \ldots, X_d > x_d) = \hat{C} \{ \bar F_1 (x_1), \ldots, \bar F_d (x_d) \}, 
\end{equation*}
where $\bar F_1, \dots, \bar F_d$ are the marginal survival functions and $\hat{C}$ is the survival copula of 
$\bs{X}$, related 
to   the copula of $\bs{X}$ as follows. If $\bs{U}$ is a random vector distributed as the copula $C$ of $\bs{X}$, 
$\hat{C}$ is 
the 
distribution function of $1-\bs{U}$. 

In risk management applications, the extremal behavior of  copulas is of particular 
interest, as it describes the dependence between extreme events and consequently the value 
of risk measures at high levels. Our purpose is to study the extremal 
behavior of Liouville copulas. The latter are defined as the survival copulas of Liouville 
distributions 
\citep{Fang/Kotz/Ng:1990, Gupta/Richards:1987, Sivazlian:1981}, 
i.e., distributions of random vectors of the 
form $R  \bs{D}_{\bs{\alpha}}$, where $R$ is a strictly positive random variable independent of the Dirichlet 
random vector $\bs{D}_{\bs{\alpha}} = (D_{1}, \dots, D_{d})$ with parameter vector $\bs{\alpha} = 
(\alpha_1, \dots, \alpha_d)$. Liouville copulas were proposed by  McNeil and Ne\v{s}lehov\'a \cite{McNeil/Neslehova:2010} 
in   order to extend the widely used class of Archimedean copulas and create dependence structures that are not necessarily 
exchangeable.  The latter 
property means that for any $u_1, \dots, u_d \in [0, 1]$ and any permutation $\pi$ of the integers $1, \dots, d$, $C(u_1, 
\dots, 
u_d)   = C(u_{\pi(1)}, \dots, u_{\pi(d)})$. When $\bs{\alpha} = \bs{1}_d \equiv(1, \dots, 1)$, 
$\bs{D}_{\bs{\alpha}}=\bs{D}_{\bs{1}_d}$ is uniformly distributed on the unit simplex 
\begin{equation}\label{eq:Sd}
\mathbb{S}_d = \{ \bs{x} \in [0, 1]^d : x_1 + \cdots + x_d = 1\}.
\end{equation}
In this special case, one recovers Archimedean copulas. Indeed, according to \cite{McNeil/Neslehova:2009}, the latter are 
the   survival copulas of random vectors  $R  \bs{D}_{\bs{1}_d}$, where $R$ is a strictly positive random variable 
independent of $\bs{D}_{\bs{1}_d}$. When $\bs{\alpha} \neq \bs{1}_d$, the survival copula of $R 
\bs{D}_{\bs{\alpha}}$ is not Archimedean anymore. It is also no longer exchangeable, unless $\alpha_1 = \dots = 
\alpha_d$.

In this article, we determine the extremal attractor of a Liouville copula and of its survival counterpart. As a by-product, we 
also obtain the lower and upper tail dependence coefficients of Liouville copulas that quantify the strength of 
dependence at extreme levels \cite{Joe:1993}. These results are complementary to \cite{Hua:2016}, where the upper tail order 
functions of a Liouville copula and its density are derived when $\alpha_1 = \dots = \alpha_d$, and to 
\cite{Hashorva:2015}, where the extremal attractor of $R  \bs{D}_{\bs{\alpha}}$ is derived when $R$ is light-tailed.
The extremal attractors of Liouville copulas are interesting in their own right. Because non-exchangeability of Liouville copulas 
carries over to their  extremal 
limits, the latter can be used to model the dependence 
between  extreme risks in the presence of causality relationships \cite{Genest/Neslehova:2013}.
The limiting extreme-value 
models can be embedded in a single family, termed here the scaled extremal Dirichlet, whose members are new, 
non-exchangeable generalizations of the 
logistic, negative logistic, and Coles--Tawn extremal Dirichlet models given in 
\cite{Coles:1991}. We examine the  scaled extremal Dirichlet model in detail and derive its de Haan spectral representation.  The 
latter is simple and leads to feasible stochastic simulation algorithms and tractable formulas for likelihood-based inference. 

The article is organized as follows. The extremal behavior of the univariate margins of Liouville distributions is 
first studied in \Cref{sec:2}. The extremal attractors of Liouville copulas and their survival counterparts are then derived in 
\Cref{sec:3}. 
When 
$\bs{\alpha}$ is integer-valued, the results of \cite{Larsson/Neslehova:2011, McNeil/Neslehova:2010} lead to closed-form 
expressions for the limiting stable tail dependence functions, as shown in \Cref{sec:4}. \Cref{sec:5} is devoted to a detailed 
study of the scaled extremal Dirichlet model. In \Cref{sec:6}, the de Haan representation is  derived 
and 
used for stochastic simulation. Estimation is investigated in \Cref{sec:7}, where expressions for the censored likelihood and 
the 
gradient score are also given. An illustrative data analysis of river flow of the river Isar  is presented in \Cref{sec:8}, and  
the paper is concluded by a discussion in \Cref{sec:9}. Lengthy proofs are relegated to 
the 
Appendices.

In what follows, vectors in $\mathbb{R}^d$ are denoted by boldface letters, $\bs{x}=(x_1, \ldots, x_d)$; $\bs{0}_d$ and 
$\bs{1}_d$  refer to the vectors $(0, \dots 0)$ and $(1, 
\ldots, 1)$ in $\mathbb{R}^d$, respectively. Binary 
operations such as $\bs{x} + \bs{y}$ or $a\cdot \bs{x}$, $\bs{x}^a$ are understood as component-wise operations.
$\|\cdot\|$   stands for the $\ell_1$-norm, viz. $\| \bs{x}\| = |x_1| + \cdots + |x_d|$, $\indep$ for statistical independence. 
For any $x, y \in \mathbb{R}$, let $x\wedge y = \min(x, y)$   and $x \vee y = \max(x, y)$.  The Dirac delta function 
$\mathrm{I}_{ij}$ is $1$ if $i=j$ and   zero otherwise. Finally, $\mathbb{R}^d_+$ is the positive orthant $[0, \infty)^d$ and for 
any $x \in 
\mathbb{R}$, $x_+$ denotes   the positive part of $x$, $\max(0, x)$.

\section{Marginal extremal behavior}\label{sec:2}
A Liouville random vector $\bs{X} = R \bs{D}_{\bs{\alpha}}$ is a scale mixture of a Dirichlet random 
vector $\bs{D}_{\bs{\alpha}} = (D_1, \dots, D_d)$ with parameters $\bs{\alpha} = 
(\alpha_1, \dots, \alpha_d)> \bs{0}_d$. In 
what follows, $R$ is referred to as the radial variable of $\bs{X}$ and $\bar \alpha$ denotes the sum of the Dirichlet 
parameters, 
viz. $\bar \alpha = \| \bs{\alpha} \| = \alpha_1 + \cdots + \alpha_d$. Recall that $\bs{D}_{\bs{\alpha}}$  has 
the 
same distribution as $\bs{Z} / \|\bs{Z}\|$, where $Z_i \sim \mathsf{Ga}(\alpha_i, 1)$, $i=1, \dots, d$ are independent 
Gamma variables with scaling parameter $1$. The margins of $\bs{X}$ are thus scale mixtures of Beta 
distributions, 
i.e., for $i=1, \dots, d$, $X_i = R D_{i}$ with $D_{i} \sim {\mathsf{Beta}}(\alpha_i, \bar \alpha - \alpha_i)$. 

As a first step towards the extremal behavior of Liouville copulas, this section is devoted to the extreme-value properties of 
the univariate margins of the vectors $\bs{X}$ and $1/\bs{X}$, where $\bs{X}$ is a Liouville random vector with 
parameters $\bs{\alpha}$ and   a strictly positive radial part $R$, i.e., such that $\Pr(R \le 0) =0$. To this end, recall that a 
univariate random 
variable $X$   with distribution function $F$ is in the maximum domain of attraction of a non-degenerate distribution $F_0$, 
denoted $F\in   \mathcal{M}(F_0)$ or $X \in \mathcal{M}(F_0)$, if and only if there exist sequences of reals  $(a_n)$ and $(b_n)$ 
with $a_n 
>0$, such that, for any $x \in \mathbb{R}$, 
\begin{align*}
\lim_{n\to\infty} F^n(a_nx + b_n) = F_0(x).                                           \end{align*}
By the Fisher--Tippett Theorem, $F_0$ must be, up to location and scale, either the Fr\'echet ($\Phi_\rho$), the Gumbel 
($\Lambda$) or the Weibull distribution ($\Psi_\rho$) with parameter $\rho > 0$. 
Further recall that a measurable function $f:\R_{+} \to \R_{+}$ is called regularly varying with index 
$\rho\in(-\infty, \infty)$, denoted $f \in \mathcal{R}_{\rho}$, if for any $x>0$, $f(tx) / f(t) \to x^\rho$ as $t \to \infty$. 
If $\rho=0$, $f$ is called slowly varying. For more details and conditions for $F\in \mathcal{M}(F_0)$, see, e.g., 
\cite{Embrechts/Kluppelberg/Mikosch:1997, Resnick:1987}.  

Because the univariate margins of $\bs{X}$ are scale mixtures of Beta distributions, their extremal behavior, 
detailed in \Cref{prop:1}, follows directly from Theorems 4.1, 4.4. and 4.5 in \cite{Hashorva/Pakes:2010}.
\begin{proposition}  
\label{prop:1}
Let $\bs{X} = R \bs{D}_{\bs{\alpha}}$ be a Liouville random vector with parameters $\bs{\alpha} = (\alpha_1, \dots, \alpha_d)$ 
and   a strictly positive radial variable $R$, i.e., $\Pr(R \le 0) =0$. Then the following statements hold for any  $\rho >0$:
\begin{enumerate}
\item[{\rm(a)}] $R\in \mathcal{M}(\Phi_\rho)$ if and only if $X_i \in \mathcal{M}(\Phi_\rho)$ for 
all $i=1, \dots, d$.
\item[{\rm(b)}] $R\in\mathcal{M}(\Lambda)$ if and only if $X_i \in \mathcal{M}(\Lambda)$ for all 
$i=1, \dots, d$.
\item[{\rm(c)}] $R\in \mathcal{M}(\Psi_\rho)$ if and only if $X_i \in \mathcal{M}(\Psi_{\rho + \bar 
\alpha - \alpha_i})$ for 
all   $i=1, \dots, d$.
\end{enumerate}
\end{proposition}
\Cref{prop:1} implies that the univariate margins of $\bs{X}$ are all in the domain of attraction of 
the same 
distribution if the latter 
is   Gumbel or Fr\'echet. This is not the case when $R$ is in the Weibull domain of attraction. 
Note also that there are 
cases   not covered by \Cref{prop:1}, in which the univariate margins $X_i$ are in the Weibull 
domain while 
$R$ is not in the 
domain 
of   attraction of any extreme-value distribution.  For example, when $d=2$, $\bs{\alpha} = (1, 1)$ 
and $R = 1$ almost surely, 
the margins of $\bs{X}$ are standard uniform and hence in the maximum domain of attraction of 
$\Psi_{1}$; see Example 3.3.15 in 
\cite{Embrechts/Kluppelberg/Mikosch:1997}. At the same time, $R$ is clearly neither in the Weibull, 
nor the Gumbel, nor the 
Fr\'echet domain of attraction.

In subsequent sections, we shall also need the extremal behavior of the univariate margins of 
$1/\bs{X}$. The 
proposition below shows that   the latter is  determined by the properties of $1/R$. In contrast to 
\Cref{prop:1}, however, the 
univariate margins of $1/\bs{X}$ are always   in the Fr\'echet domain. The proof may be found in 
\ref{app:A}.
\begin{proposition} \label{prop:2}
Let $\bs{X} = R \bs{D}_{\bs{\alpha}}$ be a Liouville random vector with parameters $\bs{\alpha} = 
(\alpha_1, \dots, \alpha_d)$ 
and   a strictly positive radial variable $R$ with $\Pr(R \le 0) =0$. The following statements hold 
for any $i=1, \dots, d$.
\begin{enumerate}
\item[\textrm{(a)}] If $1/R\in \mathcal{M}(\Phi_\rho)$ for $\rho\in(0, \alpha_i]$, then $1/X_i \in  
\mathcal{M}(\Phi_\rho)$.
\item[\textrm{(b)}] If $\e(1/R^{\alpha_i + \varepsilon}) < \infty$ for some~$\varepsilon > 0$, then 
$1/X_i \in 
\mathcal{M}(\Phi_{\alpha_i})$.
\end{enumerate}
\end{proposition}

\section{Extremal behavior of Liouville copulas}\label{sec:3}

In this section, we will identify the extremal behavior of a Liouville random vector $\bs{X} = R 
\bs{D}_{\bs{\alpha}}$ and of 
the   random vector $1/\bs{X}$, assuming that $\Pr(R \le 0) = 0$. As a by-product, we will obtain 
the extremal 
attractors of 
Liouville   copulas and their survival counterparts.
To this end, recall that a random vector $\bs{Y}$ with joint distribution function $H$ is in the 
maximum domain of attraction of 
a   non-degenerate distribution function $H_0$, in notation $H \in \mathcal{M}(H_0)$ or $\bs{Y} \in 
\mathcal{M}(H_0)$, iff 
there exist 
sequences of vectors  $(\bs{a}_{n})$ in $(0, \infty)^d$ and $(\bs{b}_{n})$ in $\mathbb{R}^d$ such 
that for all $\bs{x} \in 
\mathbb{R}^d$, 
\begin{align*}
\lim_{n\to \infty} H^n (\bs{a}_{n}\bs{x} + \bs{b}_{n} ) = H_0(\bs{x}).
\end{align*}
When the univariate margins $F_1,\dots, F_d$ of $H$ are continuous, $H \in \mathcal{M}(H_0)$ holds 
if and only if  
$F_i \in \mathcal{M}(F_{0i})$ for all $i=1, \dots, d$, where $F_{01},\dots, F_{0d}$ are the 
univariate margins of $H_0$,
and further if  the unique copula $C$ of $H$ is in the domain of attraction of the unique copula 
$C_0$ of $H_0$, denoted $C \in 
\mathcal{M}(C_0)$, i.e., iff for all $\bs{u} \in [0, 1]^d$, 
\begin{align*}
\lim_{n\to\infty} C^n(\bs{u}^{1/n}) = C_0(\bs{u}).
\end{align*}
In particular, the univariate margins of the max-stable distribution $H_0$ must each follow a 
generalized extreme-value distribution, and $C_0$ must be an 
extreme-value  copula. This means that for all $\bs{u} \in [0, 1]^d$, 
\begin{equation}\label{eq:evc}
C_0(\bs{u}) = \exp[-\ell\{-\log(u_1), \dots, -\log(u_{d})\}], 
\end{equation}
where $\ell : \mathbb{R}_+^d \to [0, \infty)$ is a stable tail dependence function, linked to the 
so-called exponent
measure $\nu$ viz. $\nu\{[\bs{0}_d, \bs{x})^\comp\}=\ell(1/\bs{x})$, see, e.g., \cite{Resnick:1987}.
The latter can be characterized through an angular (or spectral) probability 
measure $\sigma_d$ on $\Si_d$ given in \Cref{eq:Sd} which satisfies $\int_{\Si_d} w_i \, \d 
\sigma_d(\bs{w}) =1/d$ for all  
$i=1, \dots, d$. For   all $\bs{x} \in \mathbb{R}_+^d$, one has
\begin{align}
\ell(\bs{x}) = d\int_{\mathbb{S}_d}\max(w_1x_1, \dots, w_{d}x_{d}) \d \sigma_d(\bs{w}). 
\label{eq:stdfspecdens}
\end{align}
Because $\ell$ is homogeneous of order $1$, i.e., for any $c > 0$ and $\bs{x} \in \mathbb{R}_+^d$, 
$\ell(c\bs{x}) = 
c\ell(\bs{x})$, $C_0$ can also be expressed via the Pickands dependence function $\mathrm{A} : \Si_d 
\to [0, \infty)$ related to 
$\ell$ through $\ell(\bs{x}) = \|\bs{x}\| \mathrm{A}(\bs{x}/\|\bs{x}\|)$. Then at any $\bs{u} \in 
[0, 1]^d$, 
\begin{align*}
C_0(\bs{u}) = \exp\left[\log(u_1\dotsm u_d) \mathrm{A}\left\{\frac{\log(u_1)}{\log(u_1\dotsm 
u_d)}, \dots, \frac{\log(u_{d})}{\log(u_1\dotsm u_d)}\right\}\right].
\end{align*}
When $d=2$, it is more common to define the Pickands dependence function $\mathrm{A}: [0, 1] \to [0, 
1]$ through $\ell(x_1, x_2) 
=   (x_1 + x_2) \mathrm{A}\{x_2/(x_1 + x_2)\}$ so that, for all $u_1, u_2 \in [0, 1]$, 
\begin{align}\label{eq:evc2d}
C_0(u_1, u_2) = \exp\left[\log(u_1u_2) 
\mathrm{A}\left\{\frac{\log(u_2)}{\log(u_1u_2)}\right\}\right].
\end{align}

Now consider a Liouville vector $\bs{X} = R \bs{D}_{\bs{\alpha}}$ with a strictly positive radial 
variable.   
\Cref{thm:1} specifies when $\bs{X} \in \mathcal{M}(H_0)$ and identifies $H_0$.  
While part (a) follows from regular variation of $\bs{X}$, parts (b) and (c) are 
special cases of the results discussed in Section~2.2 in \cite{Hashorva:2015}. Details of the proof 
may be found in \ref{app:B}. 
\begin{theorem}\label{thm:1}
Let $\bs{X} = R \bs{D}_{\bs{\alpha}}$, $\bs{D}_{\bs{\alpha}} = (D_1, \dots, D_d)$, $\bs{\alpha} > 
\bs{0}_d$, and  
$\Pr(R \le 
0)   = 0$. Then the following statements hold. 
\begin{enumerate}
\item[\textrm{(a)}] If $R\in \mathcal{M}(\Phi_\rho)$ for some $\rho > 0$, then $\bs{X} \in 
\mathcal{M}(H_0)$, where $H_0$ is a 
multivariate extreme-value distribution with univariate margins  $F_{0i} = \Phi_\rho$, $i=1,\dots, 
d$, and a 
stable tail dependence function 
given, for all $\bs{x} 
\in   \mathbb{R}^d_+$, by
\begin{align*}
\ell(\bs{x}) = \frac{\Gamma(\bar \alpha +\rho)}{\Gamma(\bar \alpha)}{\mathrm E} \left[\max\left\{ 
\frac{\Gamma(\alpha_1) x_1 
D_{1}^{\rho}}{\Gamma(\alpha_1+\rho)}, \dots, \frac{\Gamma(\alpha_d)x_d 
D^{\rho}_{d}}{\Gamma(\alpha_d+\rho)} \right\}\right]. 
\end{align*}
\item[\textrm{(b)}] If $R \in \mathcal{M}(\Lambda)$, then $\bs{X} \in \mathcal{M}(H_0)$, where for 
all $\bs{x} \in \mathbb{R}^d$, 
$H_0(\bs{x}) = \prod_{i=1}^d\Lambda(x_i)$.
\item[\textrm{(c)}] If $R \in \mathcal{M}(\Psi_\rho)$ for some $\rho > 0$, then $\bs{X} \in 
\mathcal{M}(H_0)$, where for all 
$\bs{x}   \in \mathbb{R}^d$, $H_0(\bs{x}) = \prod_{i=1}^d\Psi_{\rho + \bar \alpha - \alpha_i}(x_i)$.
\end{enumerate}
\end{theorem}
The next result, also proved in \ref{app:B}, specifies the conditions under which $1/\bs{X} \in 
\mathcal{M}(H_0)$ and gives the 
form of the limiting extreme-value distribution $H_0$.
\begin{theorem}\label{thm:2}
Let $\bs{X} = R \bs{D}_{\bs{\alpha}}$, $\bs{D}_{\bs{\alpha}} = (D_1, \dots, D_d)$, $\bs{\alpha} > 
\bs{0}_d$, and assume that 
$\Pr(R \le 
0)   = 0$. Let $\alpha_{\mathrm{M}} = \max(\alpha_1, \dots, \alpha_d)$. The following cases can be 
distinguished:
\begin{enumerate}
\item[\textrm{(a)}] If $1/R \in  \mathcal{M}(\Phi_\rho)$ for $\rho \in(0, \alpha_{\mathrm{M}}]$, set 
$\mathbb{I}_1 = \{ i : 
\alpha_i   \le \rho\}$, $\mathbb{I}_2 = \{ i : \alpha_i > \rho\}$ and $\bar \alpha_2 = \sum_{i \in 
\mathbb{I}_2} \alpha_i$. Then 
 
$1/\bs{X}\in\mathcal{M}(H_0)$, where the univariate margins of $H_0$ are $F_{0i} =\Phi_{\rho\wedge 
\alpha_i}$, 
$i=1, \dots, d$, and the 
stable   tail dependence function is given, for all $\bs{x} \in \mathbb{R}^d_+$, by
\begin{align*}
\ell(\bs{x})  = \sum_{i \in \mathbb{I}_1} x_i + \frac{\Gamma(\bar \alpha -\rho)}{\Gamma(\bar 
\alpha)}{\mathrm E} \left[\max_{i 
\in 
\mathbb{I}_2}\left\{ \frac{\Gamma(\alpha_i)x_i D_{i}^{-\rho}}{\Gamma(\alpha_i-\rho)} 
\right\}\right] 
=\sum_{i \in \mathbb{I}_1} x_i + \frac{\Gamma(\bar \alpha_{2} -\rho)}{\Gamma(\bar 
\alpha_{2})}{\mathrm E} \left[\max_{i \in 
\mathbb{I}_2}\left\{ \frac{\Gamma(\alpha_i)x_i \widetilde D_{i}^{-\rho}}{\Gamma(\alpha_i-\rho)} 
\right\}\right], 
\end{align*}
where $(\widetilde D_{i}, i \in \mathbb{I}_2)$ is a Dirichlet random vector with parameters 
$(\alpha_i, i \in \mathbb{I}_2)$  
if $| \mathbb{I}_2| > 1$ and $\widetilde D_{i} \equiv 1$ if $\mathbb{I}_2 = \{i\}$.
\item[\textrm{(b)}] If $\E{1/R^\beta} < \infty$ for $\beta > \alpha_{\mathrm{M}}$, then $1/\bs{X} 
\in \mathcal{M}(H_0)$, where 
for 
all $\bs{x} \in \mathbb{R}^d$, $H_0(\bs{x}) = \prod_{i=1}^d\Phi_{\alpha_i}(x_i)$.
\end{enumerate}
\end{theorem}

\begin{remark}
\em Note that in the case of asymptotic independence between the components of 
$\bs{X}$ (\Cref{thm:1}~(b--c)) or $1/\bs{X}$ (\Cref{thm:2}~(b)), dependence between component-wise 
maxima of finitely many vectors may still be present. Refinements of asymptotic independence are 
then needed, but these considerations surpass the scope of this paper. One option would be to 
consider triangular arrays as in \cite{Huesler/Reiss:1989}; extremes of arrays of Liouville vectors 
can be obtained as a special case of extremes of arrays of weighted Dirichlet distributions 
developed in \cite{Hashorva:2008}. Another avenue worth exploring might be the limits of scaled 
sample clouds, as in \cite{Balkema/Nolde:2010} and \cite{Nolde:2014}.
\end{remark}

The stable tail dependence functions appearing in \Cref{thm:1,thm:2} will be investigated in greater 
detail in the subsequent 
sections. Before proceeding, we introduce the following terminology, emphasizing that they can in 
fact be embedded in one and 
the same parametric class.
\begin{definition}\label{def:1}
For any $\alpha >0$ and $\rho \in (-\alpha, \infty)$, let $c(\alpha, \rho) = \Gamma(\alpha + 
\rho)/\Gamma(\alpha)$ denote the 
rising factorial. For $d \ge 2$ and $\alpha_1, \dots, \alpha_d > 0$ and let $(D_1, \dots, D_d)$ 
denote a Dirichlet random vector 
with parameters $\bs{\alpha}=(\alpha_1, \dots, \alpha_d)$ and set $\bar \alpha = \alpha_1 + \cdots + 
\alpha_d$. For any 
$-\min(\alpha_1, \dots, \alpha_d) < \rho < \infty$, the \emph{scaled extremal Dirichlet} stable tail 
dependence 
function with parameters $\rho$ and $\bs{\alpha}$ is given, 
for all $\bs{x} \in \mathbb{R}_+^d$, by
\begin{align}\label{eq:ellD}
\ell^{\mathrm{D}}(\bs{x}; \rho, \bs{\alpha}) = c(\bar \alpha, \rho){\mathrm E} \left[\max\left\{ 
\frac{ x_1 D_{1}^{\rho} 
}{c(\alpha_1, \rho)}, \dots, \frac{x_d D^{\rho}_{d}}{c(\alpha_d, \rho)} \right\}\right], 
\end{align}
when $\rho \neq0$ and by $\max(x_1, \dots, x_d)$ when $\rho=0$.  For any $\rho > 0$, the 
\emph{positive scaled extremal 
Dirichlet} 
stable tail dependence function $\ell^{\mathrm{pD}}$ with parameters $\rho$ and $\bs{\alpha}$ is 
given, 
for all $\bs{x} \in \mathbb{R}_+^d$, by $\ell^{\mathrm{pD}}(\bs{x}; \rho, 
\bs{\alpha})=\ell^{\mathrm{D}}(\bs{x}; \rho, 
\bs{\alpha})$, while for any $0<\rho < \min(\alpha_1, \dots, \alpha_d)$, the \emph{negative scaled 
extremal Dirichlet} stable 
tail 
dependence function $\ell^{\mathrm{nD}}$ is given, for all $\bs{x} \in \mathbb{R}_+^d$, by 
$\ell^{\mathrm{nD}}(\bs{x}; \rho, 
\bs{\alpha})  = \ell^{\mathrm{D}}(\bs{x}; -\rho, \bs{\alpha})$.
\end{definition}

\begin{remark}\label{rem:1}\em
As will be seen in \Cref{sec:5}, distinguishing between the positive and negative scaled extremal 
Dirichlet models makes 
the discussion of their   properties slightly easier because the sign of $\rho$ impacts the shape of 
the corresponding angular 
measure.
When $\rho\to 0$, $\ell^{\mathrm{D}}(\bs{x}; \rho, \bs{\alpha})$  
becomes $\max(x_1, \dots, x_d)$, the stable tail 
dependence function corresponding to comonotonicity, while when $\rho\to \infty$, 
$\ell^{\mathrm{D}}(\bs{x}; \rho  
, \bs{\alpha})$ becomes $x_1 + \cdots + x_d$, the stable tail dependence function corresponding to 
independence.  Note also that 
$\rho \in (-\infty, \infty)$ can be allowed, with the convention that all variables whose indices 
$i$ are such that $\rho  \le 
-\alpha_i$ are independent, i.e., $\ell^{\mathrm{nD}}$ is then of the form given in 
\Cref{thm:2}~(a).
\end{remark}

From \Cref{thm:1,thm:2}, we can now easily deduce the extremal behavior of Liouville copulas and 
their survival counterparts. 
To   this end, recall that a Liouville copula $C$ is defined as the survival copula of a Liouville 
random vector $\bs{X} = R 
\bs{D}_{\bs{\alpha}}$ with $\Pr(R \le 0) = 0$. The following corollary follows directly from 
\Cref{thm:2} upon noting that $C$ 
is   also the unique copula of $1/\bs{X}$.
\begin{corollary}\label{cor:1}
Let $C$ be the unique survival copula of a Liouville random vector $\bs{X} = R \bs{D}_{\bs{\alpha}}$ 
with $\Pr(R \le 0) = 0$. 
Let   $\alpha_{\mathrm{M}} = \max(\alpha_1, \dots, \alpha_d)$ and set $\mathbb{I}_1 = \{ i : 
\alpha_i \le 
\rho\}$, $\mathbb{I}_2 = \{ i : \alpha_i > \rho\}$. Then the following statements hold.
\begin{enumerate}
\item[{\textrm(a)}]If $1/R \in  \mathcal{M}(\Phi_\rho)$ for $\rho \in(0, \alpha_{\mathrm{M}}]$ and 
$|\mathbb{I}_2| 
> 1$, then $C \in \mathcal{M}(C_0)$, 
where   $C_0$ is an extreme-value copula of the form \eqref{eq:evc} whose stable tail dependence 
function is given, for all 
$\bs{x} \in   \mathbb{R}^d_+$, by 
\begin{align*}
\ell(\bs{x}) = \sum_{i \in \mathbb{I}_1} x_i + \ell^{\mathrm{nD}}(\bs{x}_{\{2\}}; \rho, 
\bs{\alpha}_{\{2\}}), 
\end{align*}
where  and $\bs{x}_{\{2\}} = (x_i, i \in 
\mathbb{I}_2)$, $\bs{\alpha}_{\{2\}} = (\alpha_i, i \in \mathbb{I}_2)$. If  $1/R \in  
\mathcal{M}(\Phi_\rho)$ for $\rho \in(0, 
\alpha_{\mathrm{M}}]$ and $|\mathbb{I}_2|\le 1$, then $C \in \mathcal{M}(\Pi)$, where $\Pi$ is 
the   independence copula given, for all $\bs{u} \in [0, 1]^d$, by $\Pi(\bs{u}) = u_1\dotsm u_d$. 
\item[\textrm{(b)}]If $\E{1/R^\beta} < \infty$ for $\beta > \alpha_{\mathrm{M}}$, then $C \in 
\mathcal{M}(\Pi)$.
\end{enumerate}
\end{corollary} 

\begin{remark} \em Observe that \Cref{cor:1} (a) in particular implies that when $d=2$, $\alpha_1 < 
\alpha_2$ and 
$1/R \in  \mathcal{M}(\Phi_\rho)$ for $\alpha_1 \le \rho < \alpha_2$, $C \in \mathcal{M}(\Pi)$. Also 
note that when  $\alpha_1 = 
\dots = \alpha_d \equiv \alpha$ and $1/R \in  \mathcal{M}(\Phi_\rho)$ for 
$\rho \in(0, \alpha)$, the result in \Cref{cor:1} (a) can be derived from formula~$(5)$ in 
Proposition~3 in \cite{Hua:2016}
by relating the tail order function to the stable tail dependence function when the tail order 
equals $1$. 
\end{remark}

The survival counterpart $\hat{C}$ of a Liouville copula $C$ is given as the distribution function 
of $1-\bs{U}$, 
where $\bs{U}$ 
is   a random vector distributed as $C$. As $C$ is the unique survival copula of $\bs{X}$, $\hat{C}$ 
is the unique 
copula of 
$\bs{X}$.   The following result thus follows directly from \Cref{thm:1}.
\begin{corollary}\label{cor:2}
Let $\hat{C}$ be the unique copula of a Liouville random vector $\bs{X} = R \bs{D}_{\bs{\alpha}}$ 
with $\Pr(R \le 
0) = 0$. Then 
the   following statements hold.
\begin{enumerate}
\item[\textrm{(a)}]If $R \in  \mathcal{M}(\Phi_\rho)$ for $\rho >0$, then $\hat{C} \in 
\mathcal{M}(C_0)$, where $C_0$ 
is an 
extreme-value copula of the form \eqref{eq:evc} with the positive scaled extremal Dirichlet stable 
tail dependence function 
given, 
for all 
$\bs{x} \in   \mathbb{R}^d_+$, by 
$\ell^{\mathrm{pD}}(\bs{x};\rho, \bs{\alpha})$.
\item[\textrm{(b)}] If $R \in \mathcal{M}(\Lambda)$ or $R \in \mathcal{M}(\Psi_\rho)$ with $\rho > 
0$, then $\hat{C} 
\in 
\mathcal{M}(\Pi)$, where $\Pi$ is the independence copula.
\end{enumerate}
\end{corollary}

\section{The case of integer-valued Dirichlet parameters} \label{sec:4}
When $\bs{\alpha}$ is integer-valued, Liouville distributions are particularly tractable because 
their survival 
function   is explicit. In this section, we will use this fact to derive closed-form expressions for 
the positive and negative 
scaled   extremal Dirichlet stable tail dependence functions. To this end, first recall the notion 
of the Williamson transform. 
The latter   is related to Weyl's fractional integral transform and was used to characterize 
$d$-monotone functions in 
\cite{Williamson:1956};   it was adapted to non-negative random variables in 
\cite{McNeil/Neslehova:2009}.
\begin{definition} Let $X$ be a non-negative random variable with distribution function $F$, and let 
$k \ge 1$ be an arbitrary 
integer. The 
Williamson $k$-transform of $X$ is given, for all $x >0$, by
\begin{align*}
\mathscr{W}_k F(x) = \int_x^\infty \left ( 1- \frac{x}{r} \right)^{k-1} \d F(r) = \e 
\left(1-\frac{x}{X}\right)_+^{k-1}.
\end{align*}
\end{definition}
For any $k \ge 1$, the distribution of a positive random variable $X$ is uniquely 
determined by its Williamson $k$-transform, the formula for the inverse transform 
being explicit  \cite{McNeil/Neslehova:2009, Williamson:1956}. If $\psi = \mathscr{W}_k 
F$, then, for all $x>0$, 
\begin{align*}
F(x) = \mathscr{W}_k^{-1} \psi(x) = 1- \sum_{j=0}^{k-2} \frac{(-1)^j x^j \psi^{(j)}(x)}{j!} -  
\frac{(-1)^{k-1} x^{k-1} 
\psi^{(k-1)}_+(x)}{(k-1)!}, 
\end{align*}
where for $j = 1, \dots, k-2$, $\psi^{(j)}$ is the $j$th derivative of $\psi$ and $\psi^{(k-1)}_+$  
is the right-hand 
derivative 
of~$\psi^{(k-2)}$. These derivatives exist because a Williamson $k$-transform $\psi$ is necessarily 
$k$-monotone 
\cite{Williamson:1956}. This means that $\psi$ is differentiable up to order $k-2$ on $(0, \infty)$ 
with derivatives satisfying 
$(-1)^j \psi^{(j)} \ge 0$ for $j=0, \dots, k-2$ and such that $(-1)^{k-2} \psi^{(k-2)}$ is 
non-increasing and convex on 
$(0, \infty)$. Moreover, $\psi(x) \to 0$ as $x \to \infty$ and if $F(0)=0$, $\psi(x) \to 1$ and $x 
\to 0$.

Now let $C$ be a Liouville copula corresponding to a Liouville random vector $\bs{X} =R 
\bs{D}_{\bs{\alpha}}$ with 
integer-valued   parameters $\bs{\alpha}=(\alpha_1, \dots, \alpha_d)$ and a strictly positive radial 
part $R$, i.e., $\Pr(R \le 
0) 
= 0$. Let $\psi$   be the Williamson $\bar \alpha$-transform of $R$ and set 
$\mathbb{I}_{\bs{\alpha}}=\{0, \dots, \alpha_1 -1\} 
\times \dots \times    \{0, \dots, \alpha_d -1\}$.  By Theorem~2 in  \cite{McNeil/Neslehova:2010}, 
one then has, for all 
$\bs{x} \in   \mathbb{R}_+^d$, 
\begin{equation}\label{eq:1}
\Pr(\bs{X} > \bs{x}) = \bar H(\bs{x})=\sum_{\mathclap{(j_1, \dots, j_d) \in 
\mathbb{I}_{\bs{\alpha}}}}\;\;(-1)^{j_1 + \cdots + 
j_d}   \frac{\psi^{(j_1 + \cdots + j_d)}(x_1 + \dotsm +x_d)}{j_1! \dotsm j_d!} \prod_{i=1}^d 
x_i^{j_i}.
\end{equation}
In particular, the margins of $\bs{X}$ have survival functions satisfying, for all $x >0$ and $i=1, 
\dots, d$, 
\begin{equation}\label{eq:2}
\Pr(X_i > x) = \bar H_i (x)= \sum_{j = 0}^{\alpha_i -1} \frac{(-1)^j x^j\psi^{(j)}(x) }{j!} = 1- 
\mathscr{W}^{-1}_{\alpha_i} 
\psi(x).
\end{equation}
By Sklar's Theorem for survival functions, the Liouville copula $C$ is given, for all $\bs{u} \in 
[0, 1]^d$, by
\begin{align*}
C(\bs{u}) = \bar H\{ \bar H_1^{-1}(u_1), \dots, \bar H_d^{-1}(u_d)\}.                                
\end{align*}
Although this formula is not explicit, it is clear from Equations \eqref{eq:1} and \eqref{eq:2} that 
$C$ depends on the 
distribution of $\bs{X}$ only through the Williamson $\bar \alpha$-transform $\psi$ of $R$ and the 
Dirichlet 
parameters~$\bs{\alpha}$.  For this reason, we shall denote the Liouville copula in this section by 
$C_{\psi, \bs{\alpha}}$ and 
refer to $\psi$ as its generator, reiterating that $\psi$ must be an $\bar \alpha$-monotone function 
satisfying $\psi(1) =0$ and 
$\psi(x) \to 0$ as $x \to \infty$. When $\bs{\alpha} = \bs{1}_d$, $C_{\psi, \bs{1}}$ is the 
Archimedean copula with 
generator $\psi$, given, for all $\bs{u} \in [0, 1]^d$ by
$
C_{\psi, \bs{1}}(\bs{u}) = \psi\{ \psi^{-1}(u_1) + \cdots + \psi^{-1}(u_d)\}.                        
$
Because the relationship between $\psi$ and $R$ is one-to-one \citep[Proposition 
3.1]{McNeil/Neslehova:2009}, we will refer to 
$R$   as the radial distribution corresponding to $\psi$. 

Now suppose that $1/R \in \mathcal{M}(\Phi_\rho)$ with $\rho \in (0, 1)$. By Theorem~2 in 
\cite{Larsson/Neslehova:2011}, this 
condition is equivalent to $1-\psi(1/\cdot) \in \mathcal{R}_{-\rho}$. It further follows from 
\Cref{cor:1}~(a) that 
$C_{\psi, \bs{\alpha}} \in \mathcal{M}(C_0)$ where $C_0$ is an extreme-value copula with the 
negative scaled extremal Dirichlet 
stable tail dependence function $\ell^{\mathrm{nD}}(\cdot; \rho, \bs{\alpha})$. This is because 
$\rho < 1 \le \min(\alpha_1, \dots, \alpha_d)$ so that 
$\mathbb{I}_1=\emptyset$ in \Cref{cor:1} (a). \Cref{eq:1} and the results of 
\cite{Larsson/Neslehova:2011} can now be used to 
derive the following explicit expression for $\ell^{\mathrm{nD}}$, as detailed in \ref{app:C}.
\begin{proposition}\label{prop:3}
Let $C_{\psi, \bs{\alpha}}$ be a Liouville copula with integer-valued parameters 
$\bs{\alpha}=(\alpha_1, \dots, \alpha_d)$ 
and generator $\psi$. If $1- \psi(1/\cdot) \in \mathcal{R}_{-\rho} $ for some $\rho \in 
(0, 1)$, then  $C_{\psi, \bs{\alpha}} \in \mathcal{M}(C_0)$, where $C_0$ is an extreme-value copula 
with scaled negative 
extremal 
Dirichlet stable tail dependence function $\ell^{\mathrm{nD}}$ as given in \Cref{def:1}. 
Furthermore, for all 
$\bs{x} \in \mathbb{R}^d_+$, 
\begin{align*}
\ell^{\mathrm{nD}}(\bs{x};\rho, \bs{\alpha}) = \Gamma(1-\rho)\left[
\sum_{j=1}^d\left\{\frac{x_j}{c(\alpha_j, -\rho)} \right\}^{1/\rho}
\right]^{\, \rho} 
\left( 1- \rho\;\;\sum_{\mathclap{\substack{(j_1, \dots, j_d ) \in 
\mathbb{I}_{\bs{\alpha}} \\ (j_1, \dots, j_d) \neq (0, \dots, 
0)}}}\quad \frac{\Gamma(j_1 + \cdots + j_d - \rho)}{\Gamma(1-\rho)} \prod_{i=1}^d 
\frac{1}{\Gamma(j_i+1)}\left[ 
\frac{\left\{\frac{x_i}{c(\alpha_i, -\rho)}\right\}^{1/\rho}}{
\sum_{k=1}^d \left\{\frac{x_k}{c(\alpha_k, -\rho)}\right\}^{1/\rho}}
\right]^{j_i} \right).
\end{align*}
\end{proposition}
When $\bs{\alpha} = \bs{1}_d$, the index set $\mathbb{I}_{\bs{\alpha}}$ reduces to the singleton 
$\{\bs{0}\}$, and 
the 
expression for $\ell^{\mathrm{nD}}$ given in \Cref{prop:3} simplifies, for all $\bs{x} \in 
\mathbb{R}_+^d$, to the  
stable tail dependence function of the Gumbel--Hougaard copula, viz. 
\begin{align*}
\ell^{\mathrm{nD}}(\bs{x}; \rho, \bs{1}_d) = \bigl(x_1^{1/\rho} + \cdots + 
x_d^{1/\rho}\bigr)^{\rho}.                             
\end{align*}
The Liouville copula $C_{\psi, \bs{1}_d}$, which is the Archimedean copula with generator $\psi$, is 
thus indeed in the domain of 
attraction of the Gumbel--Hougaard copula with parameter $1/\rho$, as shown, e.g., in 
\cite{Charpentier/Segers:2009, Larsson/Neslehova:2011}.
\begin{remark} \em
When $\bs{\alpha} = \bs{1}_d$ and $1-\psi(1/\cdot) \in \mathcal{R}_{-1}$, it is shown in 
Proposition~2 of 
\cite{Larsson/Neslehova:2011}  that $C_{\psi, \bs{1}}$ is in the domain of attraction of the 
independence copula. However, when 
$\bs{\alpha}$ is integer-valued but such that $\max(\alpha_1, \dots, \alpha_d) > 1$, regular 
variation of $1-\psi(1/\cdot)$ 
does 
not suffice to characterize those cases in \Cref{cor:1} that are not covered by \Cref{prop:3}.  This 
is because by 
Theorem~2 of \cite{Larsson/Neslehova:2011}, $1/R \in \mathcal{M}(\Phi_\rho)$ for $\rho \ge 1$, $1/R 
\in \mathcal{M}(\Lambda)$ 
and 
$1/R \in \mathcal{M}(\Psi_\rho)$ for $\rho >0$ all imply that $1-\psi(1/\cdot) \in 
\mathcal{R}_{-1}$. At the same time, by 
\Cref{cor:1}, $C_{\psi, \bs{\alpha}}\in\mathcal{M}(\Pi)$ clearly does not hold in all these cases.
\end{remark} 
Next, let $\hat{C}_{\psi, \bs{\alpha}}$ be the survival copula of a Liouville copula $C_{\psi, 
\bs{\alpha}}$, 
i.e., the 
distribution function of $1-\bs{U}$, where $\bs{U}$ is a random vector with distribution function 
$C_{\psi, \bs{\alpha}}$. The 
results of \cite{Larsson/Neslehova:2011} can again be used to restate the conditions under which 
$\hat{C}_{\psi, 
\bs{\alpha}} 
\in 
\mathcal{M}(C_0)$ in terms of $\psi$ and to give an explicit expression for the stable tail 
dependence function of $C_0$. 

\begin{proposition}\label{prop:4}
Let $\hat{C}_{\psi, \bs{\alpha}}$ be the survival copula of a Liouville copula $C_{\psi, 
\bs{\alpha}}$ with 
integer-valued 
parameters 
$\bs{\alpha}$ and a generator $\psi$. Then the following statements hold. 
\begin{enumerate}
\item[\textrm{(a)}] If $\psi\in\mathcal{R}_{-\rho}$ for some $\rho > 0$, then $\hat{C}_{\psi, 
\bs{\alpha}} \in 
\mathcal{M}(C_0)$, 
where 
$C_0$ has a positive scaled extremal Dirichlet stable tail dependence function $\ell^{\mathrm{pD}}$ 
as given in \Cref{def:1}. 
The 
latter can be expressed, for all $\bs{x} \in \mathbb{R}_+^d$, as
\begin{multline*}
\ell^{\mathrm{pD}}(\bs{x};\rho, \bs{\alpha}) = \frac{\Gamma(1+\rho)}{\Gamma(\rho)}\sum_{k=1}^d 
\sum_{1 \le i_1 < \dots <  i_k\le 
d} 
\!\!\!(-1)^{k+1} 
\left[
\left\{\sum_{j=1}^k \frac{x_{i_j}}{c(\alpha_{i_j}, \rho)}\right\}^{-1/\rho}
\right]^{-\rho} \times \\
\sum_{(j_1, \dots, j_k) \in \mathbb{I}_{(\alpha_{i_1}, \dots, \alpha_{i_k})}} 
\!\!\!\frac{\Gamma(j_1 + \cdots + j_k + \rho)}{ j_1! \dotsm 
j_k!} \prod_{m=1}^k 
\left[\frac{\left\{\frac{x_{i_m} }{c(\alpha_{i_m},\,
\rho)}\right\}^{-1/\rho}}{
\sum_{j=1}^k \left\{\frac{x_{i_j}}{c(\alpha_{i_j},\, \rho)}\right\}^{-1/\rho}
}\right]^{j_m}.
\end{multline*}
\item[\textrm{(b)}] If $\psi \in \mathcal{M}(\Lambda)$ or $\psi \in \mathcal{M}(\Psi_\rho)$ for some 
$\rho > 0$, $\bar 
C_{\psi, \bs{\alpha}} \in \mathcal{M}(\Pi)$, where $\Pi$ is the independence copula.
\end{enumerate}
\end{proposition}
When $\bs{\alpha} = \bs{1}_d$, the expression for $\ell^{\mathrm{pD}}$ in part (a) of \Cref{prop:4} 
simplifies, for all $\bs{x} 
\in 
\mathbb{R}^d_{+}$, to
\begin{align*}
\ell^{\mathrm{pD}}(\bs{x};\rho, \bs{1}_d) = \sum_{\substack{A \subseteq \{1, \dots, d\}, \; A \neq 
\emptyset}} (-1)^{|A| + 
1}\left( 
\sum_{i \in A} x_i^{-1/\rho} \right)^{-\rho}, 
\end{align*}
which is the stable tail dependence function of the Galambos copula \cite{Joe:1990b}. When $\psi \in 
\mathcal{R}_{-\rho}$ for 
some 
$\rho >0$, $\hat{C}_{\psi, \bs{1}_d}$ is thus indeed in the domain of attraction of the Galambos 
copula, as shown, 
e.g., in  
\cite{Larsson/Neslehova:2011}.

\section{Properties of the scaled extremal Dirichlet models}\label{sec:5}
In this section, the scaled extremal Dirichlet model with stable tail dependence function given in 
\Cref{def:1} is investigated in greater detail. In \Cref{sec:5a} we derive formulas for the 
so-called angular density and  
relate 
the positive and negative scaled extremal Dirichlet models to classical classes of stable tail 
dependence functions. In 
\Cref{sec:5b} we focus on the bivariate case and derive explicit expressions for the stable tail 
dependence functions and, as 
a by-product, obtain formulas for the  tail dependence coefficients of Liouville copulas.

\begin{figure}[t!]
\centering 
\includegraphics[width=0.7\textwidth]{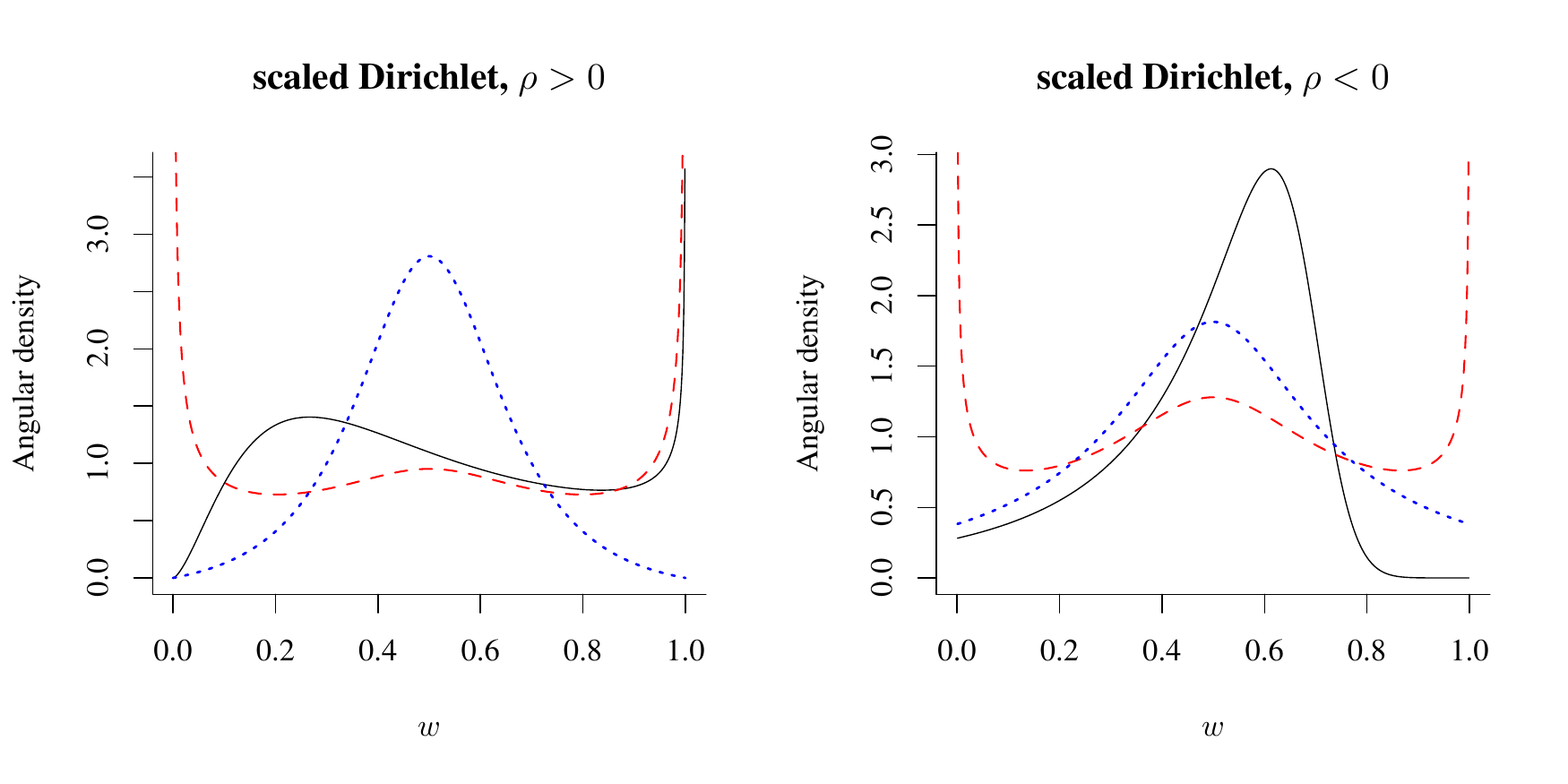}
\caption{Angular density of the scaled extremal Dirichlet model. Left panel:  $\rho=4/5$ and 
$\bs{\alpha}=(2, 1/2)$ 
(black full), $\rho=1/4$ and $\bs{\alpha}=(1/10, 1/10)$ (red dashed), $\rho=1/4$ and 
$\bs{\alpha}=(1/2, 1/2)$ (blue dotted). 
Right panel:  $\rho=-1/4$ and $\bs{\alpha}=(2, 1/2)$ (black full), $\bs{\alpha}=(2/5, 2/5)$ (red 
dashed) and $\bs{\alpha}=(1/2, 
1/2)$ (blue dotted). }
\label{fig:angular2d}
\end{figure}

\subsection{Angular density} \label{sec:5a}
The first property worth noting is that the positive and negative scaled extremal Dirichlet models 
are closed  under 
marginalization. 
Indeed, letting $x_i \to 0$ for some arbitrary $1 \le i \le d$, we can 
easily derive from \Cref{lem:B1} that for any $\bs{\alpha} > \bs{0}_d$, $\rho >0$, and any $\bs{x} 
\in \mathbb{R}_+^d$, 
$
\ell^{\mathrm{pD}}(\bs{x}; \rho, \bs{\alpha}) \to \ell^{\mathrm{pD}}(\bs{x}_{-i}; \rho, 
\bs{\alpha}_{-i})
$ as $x_i \to 0$, 
where for any $\bs{y} \in \mathbb{R}^d$, $\bs{y}_{-i}$ denotes the vector $(y_1, \dots, y_{i-1}, 
y_{i+1}, \dots, y_d)$. 
Similarly, 
for any $1\le i \le d$, $\bs{\alpha} > \bs{0}_d$, $0< \rho < \min(\alpha_1, \dots, \alpha_d)$, and 
any $\bs{x} \in 
\mathbb{R}_+^d$, 
$ \ell^{\mathrm{nD}}(\bs{x}; \rho, \bs{\alpha}) \to \ell^{\mathrm{nD}}(\bs{x}_{-i}; \rho, 
\bs{\alpha}_{-i})$ as $x_i \to 0$.

Since none of the scaled extremal Dirichlet models places mass on the vertices or facets of the 
simplex $\mathbb{S}_d$  when 
$\rho \neq 0$, the 
density of the angular measure $\sigma_d$ completely characterizes the stable tail dependence 
function and hence also the 
associated extreme-value copula. This so-called angular density of the scaled extremal Dirichlet 
models 
is 
given below and derived in \ref{app:D}.
\begin{proposition}\label{prop:5}   Let $d \ge 2$ and set $\alpha_1, \dots, \alpha_d > 0$ and $\bar 
\alpha = \alpha_1 + \cdots + 
\alpha_d$. For any  $\rho > - \min(\alpha_1, \dots, \alpha_d)$, let also $\bs{c}(\bs{\alpha}, \rho) 
= (c(\alpha_1, \rho), \dots, 
c(\alpha_d, \rho))$, where $c(\alpha, \rho)$ is as in \Cref{def:1}. Then for any $-\min(\alpha_1, 
\dots, \alpha_d) < \rho < 
\infty$, 
$\rho \neq 0$, the angular density of the scaled extremal Dirichlet model with parameters $\rho >0$ 
and 
$\bs{\alpha}$ 
is given, for all $\bs{w} \in \mathbb{S}_d$, by
\begin{align*}
      h^{\mathrm{D}}(\bs{w};\rho, 
\bs{\alpha})=\frac{\Gamma(\bar{\alpha}+\rho)}{d|\rho|^{d-1}\prod_{i=1}^d\Gamma(\alpha_i)}
\left[\sum_{j=1}^d\left\{c(\bs{\alpha}_j, 
\rho)w_j\right\}^{1/\rho}\right]^{-\rho-\bar{\alpha}}\prod_{i=1}^{d} 
\{c(\alpha_i, \rho)\}^{\alpha_i/\rho}w_i^{\alpha_i/\rho-1}.
\end{align*}
The angular density of the positive scaled extremal Dirichlet model with parameters $\rho > 0$ and 
$\bs{\alpha}$ is given, for 
all
$\bs{w} \in \mathbb{S}_d$, by $h^{\mathrm{pD}} (\bs{w};\rho, \bs{\alpha})= 
h^{\mathrm{D}}(\bs{w};\rho, \bs{\alpha})$, while the 
angular density of the negative scaled extremal Dirichlet model with parameters $0< \rho < 
\min(\alpha_1, \dots, \alpha_d)$ and 
$\bs{\alpha}$ is given, for all $\bs{w} \in \mathbb{R}_+^d$, by $h^{\mathrm{nD}} (\bs{w};\rho, 
\bs{\alpha})= 
h^{\mathrm{D}}(\bs{w};-\rho, \bs{\alpha})$.
\end{proposition}
From \Cref{prop:5}, it is easily seen that when $\bs{\alpha} = \bs{1}_d$, the angular density 
$h^{\mathrm{pD}}$ 
reduces, for any $\rho > 0$ and $\bs{w} \in \mathbb{S}_d$, to the angular density of the symmetric 
negative logistic model; 
see, e.g., Section~4.2 in \cite{Coles:1991}. In general, the angular density $h^{\mathrm{pD}}$ is 
not symmetric unless 
$\bs{\alpha}=\alpha\bs{1}_d$. 

The positive scaled Dirichlet model can thus be viewed as a new asymmetric generalization of the 
negative 
logistic model which does not place any mass on the vertices or facets of $\mathbb{S}_d$, unless at 
independence or 
comonotonicity, i.e., when $\rho \to \infty$ and $\rho \to 0$, respectively. Furthermore, 
$h^{\mathrm{pD}}$ can also be 
interpreted 
as a generalization of the Coles--Tawn extremal Dirichlet model. Indeed, $h^{\mathrm{pD}}(\bs{x};1, 
\bs{\alpha})$ is precisely 
the 
angular density of the latter model given, e.g., in Equation (3.6) in \cite{Coles:1991}. Similarly, 
the negative scaled 
extremal Dirichlet model is a new asymmetric generalization of Gumbel's logistic 
model  \cite{Gumbel:1960}. Indeed, when $\bs{\alpha} =\bs{1}_d$, $h^{\mathrm{nD}}$ simplifies to the 
logistic angular density, 
given, e.g., on p. 381 in \cite{Coles:1991}. 

\begin{figure}[t!]
\centering 
\includegraphics[width=0.75\textwidth]{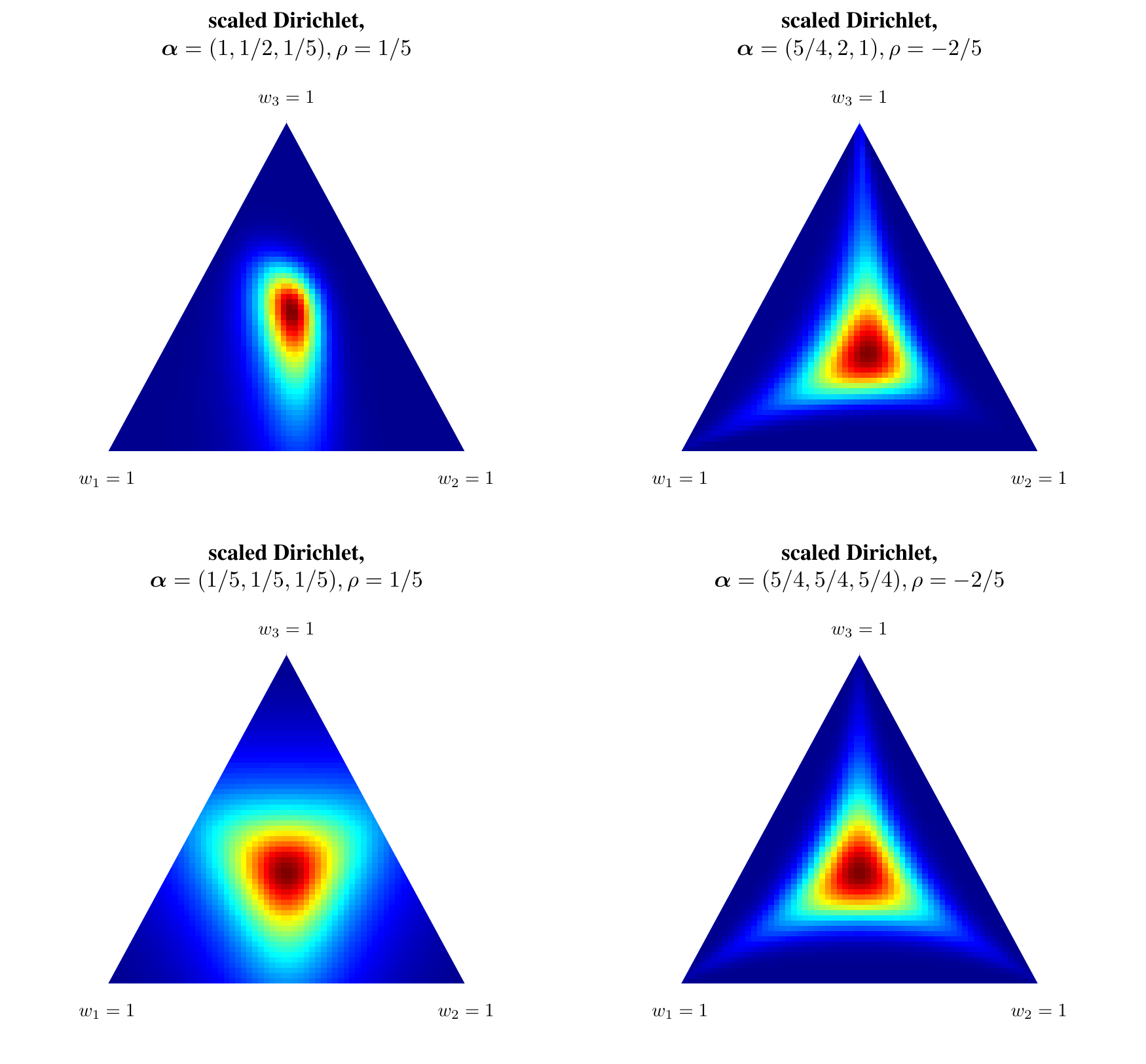}
\caption{Angular density of the scaled extremal Dirichlet model with $\bs{\alpha}=(1, 1/2, 1/5), 
\rho=1/5$ (top left) 
and 
$\bs{\alpha}=(1/5, 1/5, 1/5), \rho=1/5$ (bottom left), 
$\bs{\alpha}=(5/4, 2, 1)$, $\rho=-2/5$ (top right) and $\bs{\alpha}=(5/4, 5/4, 5/4), \rho=-2/5$ 
(bottom right). The colors 
correspond 
to log density values and range from red (high density) to 
blue (low density).}
\label{fig:angular3d}
\end{figure}

Figures \ref{fig:angular2d} and \ref{fig:angular3d} illustrate the various shapes of 
$h^{\mathrm{pD}}$ and $h^{\mathrm{nD}}$ 
that 
obtain through various choices of $\bs{\alpha}$ and $\rho$. The asymmetry when $\bs{\alpha} \neq 
\alpha \bs{1}_d$ is clearly 
apparent. For the same value of $\rho$, the shapes of the angular density can be quite different 
depending on $\bs{\alpha}$. In view of the aforementioned closure of both the positive and negative 
scaled extremal Dirichlet 
models under marginalization, this means that these models are able to capture strong dependence in 
some pairs of variables 
(represented by a mode close to $1/2$ of the angular density) and at the same time weak dependence 
in others pairs (represented 
by a bathtub shape).

\subsection{The bivariate case} \label{sec:5b}

When $d=2$, the stable tail dependence functions of the positive and negative scaled extremal 
Dirichlet models have a 
closed-form 
expression  in terms of the incomplete beta function given, for any $t \in (0, 1)$ and $\alpha_1, 
\alpha_2 > 0$, by
\begin{align*}
\Be(t;\alpha_1, \alpha_2) = \int_0^t x^{\alpha_1 - 1} (1-x)^{\alpha_2 -1} \d x.                      
\end{align*}
When $t=1$, this integral is the beta function, viz. $\Be(\alpha_1, \alpha_2) = \Gamma(\alpha_1) 
\Gamma(\alpha_2)/\Gamma(\alpha_1+\alpha_2)$. A direct calculation yields 
the corresponding Pickands dependence function, for any $t \in [0, 1]$, 
$\mathrm{A}^{\mathrm{pD}}(t; 
\rho, \alpha_1, \alpha_2)   = \ell^{\mathrm{pD}}(1-t, t;\rho, \alpha_1, \alpha_2)$, i.e., 

\begin{multline*}
\mathrm{A}^{\mathrm{pD}}(t; \rho, \alpha_1, \alpha_2)   =  \frac{(1-t)}{\Be(\alpha_2, \alpha_1 + 
\rho)} \Be\left[ 
\frac{\left\{c(\alpha_2, \rho)(1-t)\right\}^{1/\rho}}{\left\{c(\alpha_2, 
\rho)(1-t)\right\}^{1/\rho}+\left\{c(\alpha_1, 
\rho)t\right\}^{1/\rho}}; 
\alpha_2, \alpha_1 + \rho\right] \\
+ \frac{t}{\Be(\alpha_1, \alpha_2+\rho)} \Be\left[ \frac{\left\{c(\alpha_1, 
\rho)t\right\}^{1/\rho}}{\left\{c(\alpha_2, 
\rho)(1-t)\right\}^{1/\rho}+\left\{c(\alpha_1, \rho)t\right\}^{1/\rho}}; \alpha_1, 
\alpha_2+\rho\right].
\end{multline*}
When $\alpha_1 =\alpha_2 = 1$, $\mathrm{A}^{\mathrm{pD}}$ becomes the Pickands dependence function 
of the Galambos copula, viz.
$\mathrm{A}^{\mathrm{pD}}(t;\rho,1,1) = 1- \bigl\{ t^{-1/\rho} + (1-t)^{-1/\rho}\bigr\}^{-\rho}$, as 
expected given that the 
positive scaled extremal 
Dirichlet 
model becomes the symmetric negative logistic model in this case.

\begin{figure}[t!]
\centering 
\includegraphics[width=0.8\textwidth]{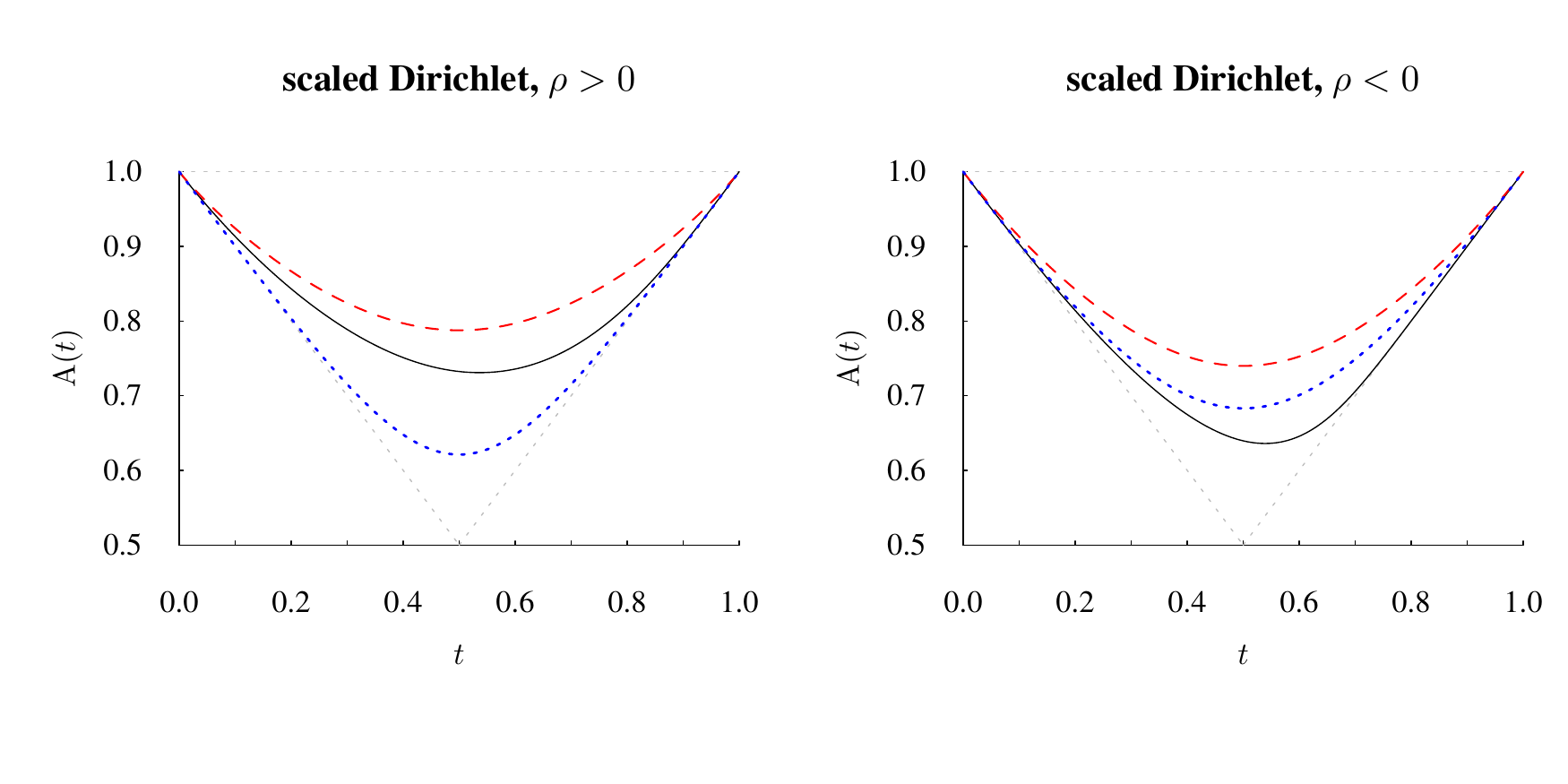}
\caption{Pickands dependence function of the 
scaled extremal Dirichlet model. Left panel:  $\rho=4/5$ and $\bs{\alpha}=(2, 1/2)$ 
(black full), $\rho=1/4$ and $\bs{\alpha}=(1/10, 1/10)$ (red dashed), $\rho=1/4$ and 
$\bs{\alpha}=(1/2, 1/2)$ (blue dotted). 
Right panel:  $\rho=-1/4$ and $\bs{\alpha}=(2, 1/2)$ (black full), $\bs{\alpha}=(2/5, 2/5)$ (red 
dashed) and $\bs{\alpha}=(1/2, 
1/2)$ (blue dotted).}
\label{fig:pickands}
\end{figure}

Similarly, for any $t \in [0, 1]$,
the Pickands dependence function $\mathrm{A}^{\mathrm{nD}}(t;\rho, \alpha_1, \alpha_2)   = 
\ell^{\mathrm{nD}}(1-t, t;\rho, \alpha_1, \alpha_2)$ equals
\begin{multline*}
\mathrm{A}^{\mathrm{nD}}(t;\rho, \alpha_1, \alpha_2)   = \frac{(1-t)}{\Be(\alpha_1-\rho, \alpha_2)} 
\Be\left[ 
\frac{\left\{(1-t)c(\alpha_2, - 
\rho)\right\}^{1/\rho}}{\left\{tc(\alpha_1, - \rho)\right\}^{1/\rho}+\left\{(1-t)c(\alpha_2, - 
\rho)\right\}^{1/\rho}}; \alpha_1 
-\rho, \alpha_2\right] \\
+ \frac{t}{\Be(\alpha_2-\rho, \alpha_1)} \Be\left[ \frac{\left\{tc(\alpha_1, - 
\rho)\right\}^{1/\rho}}{\left\{tc(\alpha_1, - 
\rho)\right\}^{1/\rho}+\left\{(1-t)c(\alpha_2, - \rho)\right\}^{1/\rho}}; \alpha_2-\rho, 
\alpha_1\right].
\end{multline*}
When $\alpha_1 = \alpha_2 = 1$, $\mathrm{A}^{\mathrm{nD}}$ simplifies to the stable tail dependence 
function of the Gumbel 
extreme-value 
copula, viz.
$\mathrm{A}^{\mathrm{nD}}(t;\rho,1,1) = \bigl\{t^{1/\rho} + (1-t)^{1/\rho}\bigr\}^{\rho}$.
This again confirms that the negative scaled extremal Dirichlet model becomes the symmetric logistic 
model when $\alpha_1= 
\alpha_2=1$.
The Pickands dependence functions $\mathrm{A}^{\mathrm{pD}}$ and $\mathrm{A}^{\mathrm{nD}}$ are 
illustrated in Figure 
\ref{fig:pickands}, 
for the same 
choices of parameters and the corresponding angular density shown in Figure \ref{fig:angular2d}.

The above formulas for  $\mathrm{A}^{\mathrm{pD}}$ and $\mathrm{A}^{\mathrm{nD}}$ now easily lead to 
expressions for their upper 
tail 
dependence 
coefficients. Recall that for an arbitrary bivariate copula $C$, the lower and upper tail dependence 
coefficients of 
\cite{Joe:1993} are given by
\begin{align*}
\lambda_\ell(C) = \lim_{u \to 0} \frac{C(u, u)}{u}, \quad \lambda_u(C) = 2- \lim_{u \to 1} 
\frac{C(u, u)-1}{u-1} = \lim_{u \to 
0} 
\frac{\hat{C}(u, u)}{u}, 
\end{align*}
where $\hat{C}$ is the survival copula of $C$, provided these limits exist. When $C$ is bivariate 
extreme-value 
with Pickands 
dependence function $\mathrm{A}$, it follows easily from \eqref{eq:evc2d} that $\lambda_\ell(C)=0$ 
and $\lambda_u(C) = 
2-2\mathrm{A}(1/2)$.

Now suppose that $C^{\mathrm{pD}}_{\rho, \bs{\alpha}}$ is a bivariate extreme-value copula with 
positive scaled extremal 
Dirichlet Pickands dependence function $\mathrm{A}^{\mathrm{pD}}$ and parameters $\rho > 0$ and $ 
\alpha_1, \alpha_2 > 0$. 
Then
\begin{multline}\label{eq:tailpD}
\lambda_u(C^{\mathrm{pD}}_{\rho, \bs{\alpha} })  = 2- \frac{1}{\Be(\alpha_2, \alpha_1 + \rho)} 
\Be\left\{ \frac{c(\alpha_2, 
\rho)^{1/\rho}}{c(\alpha_2, \rho)^{1/\rho}+c(\alpha_1, \rho)^{1/\rho}}; \alpha_2, \alpha_1 + 
\rho\right\} \\
- \frac{1}{\Be(\alpha_1, \alpha_2+\rho)} \Be\left\{ \frac{c(\alpha_1, \rho)^{1/\rho}}{c(\alpha_2, 
\rho)^{1/\rho}+c(\alpha_1, 
 \rho)^{1/\rho}}; \alpha_1, \alpha_2+\rho\right\}.
\end{multline}
Similarly, if $C^{\mathrm{nD}}_{\rho, \bs{\alpha}}$ is a bivariate extreme-value copula with  
negative scaled extremal 
Dirichlet Pickands dependence function $\mathrm{A}^{\mathrm{nD}}$ and parameters $\alpha_1, \alpha_2 
> 0$ and $0 < \rho < 
\min(\alpha_1, \alpha_2)$, then
\begin{multline}\label{eq:tailnD}
\lambda_u(C^{\mathrm{nD}}_{\rho, \bs{\alpha} })   = 2- \frac{1}{\Be(\alpha_1-\rho, \alpha_2)} 
\Be\left\{ \frac{c(\alpha_2, - 
\rho)^{1/\rho}}{c(\alpha_1, - \rho)^{1/\rho}+c(\alpha_2, - \rho)^{1/\rho}}; \alpha_1 -\rho, 
\alpha_2\right\} \\
+ \frac{1}{\Be(\alpha_2-\rho, \alpha_1)} \Be\left\{ \frac{c(\alpha_1, - \rho)^{1/\rho}}{c(\alpha_1, 
- \rho)^{1/\rho}+c(\alpha_2 
, - \rho)^{1/\rho}}; \alpha_2-\rho, \alpha_1\right\}.
\end{multline}
In the symmetric case $\alpha_1 = \alpha_2\equiv \alpha$, Expressions \eqref{eq:tailpD} and 
\eqref{eq:tailnD} simplify to 
\begin{align*}
\lambda_u(C^{\mathrm{pD}}_{\rho, \bs{\alpha} }) = 2-\frac{2} {\Be(\alpha, \alpha + \rho)} \Be\left( 
\frac{1}{2}; \alpha, \alpha + 
\rho\right), \quad  \lambda_u(C^{\mathrm{nD}}_{\rho, \bs{\alpha} }) = 2-\frac{2} {\Be(\alpha-\rho, 
\alpha )} \Be\left( 
\frac{1}{2}; 
\alpha-\rho, \alpha\right).
\end{align*}
Formulas \eqref{eq:tailpD} and \eqref{eq:tailnD} lead directly to expressions for the tail 
dependence coefficients of Liouville 
copulas.  This is because if $C \in \mathcal{M}(C_0)$, where $C_0$ is an extreme-value copula with 
Pickands tail dependence 
function $\mathrm{A}_0$, $\lambda_u(C)= 2-2\mathrm{A}_0(1/2)$ \citep[Proposition 
7.51]{McNeil/Frey/Embrechts:2015}. Similarly, if 
$\hat{C} \in 
\mathcal{M}(C_0^*)$, where $C^*_0$ is an extreme-value copula with Pickands tail dependence function 
$\mathrm{A}^*_0$, 
$\lambda_\ell(C)= 
2-2\mathrm{A}^*_0(1/2)$. The following corollary is thus an immediate consequence of Corollaries 
\ref{cor:1} and \ref{cor:2}.

\begin{corollary}\label{cor:3}
Suppose that $C$ is the survival copula of a Liouville random vector $R\bs{D}_{\bs{\alpha}}$ with 
parameters $\bs{\alpha} > 0$ 
and 
a radial part $R$ such that $\Pr(R \le 0) =0$. Then the following statements hold.
\begin{enumerate}
\item[\textrm{(a)}] If $R \in \mathcal{M}(\Phi_\rho)$ for some $\rho > 0$, $\lambda_\ell(C) = 
\lambda_u(C^{\mathrm{pD}}_{\rho, \bs{\alpha}})$ is given by \Cref{eq:tailpD}. 
\item[\textrm{(b)}] If $R \in \mathcal{M}(\Lambda)$ or $R \in \mathcal{M}(\Psi_\rho)$ for some $\rho 
> 0$, $\lambda_\ell(C)=0$.
\item[\textrm{(c)}] If $1/R \in \mathcal{M}(\Phi_\rho)$ for some $0 < \rho < \alpha_1 \wedge 
\alpha_2$, $\lambda_u(C) = 
\lambda_u(C^{\mathrm{nD}}_{\rho, \bs{\alpha}})$ is given by \Cref{eq:tailnD}. 
\item[\textrm{(d)}] If $1/R \in \mathcal{M}(\Phi_\rho)$ for $\rho > \alpha_1 \wedge \alpha_2$ or if  
$\e({1/R^\beta}) < \infty$ 
for 
$\beta> \alpha_1 \vee \alpha_2$, $\lambda_u(C) =0$.
\end{enumerate}
\end{corollary}

The role of the parameters $\bs{\alpha}$ and $\rho$ is best explained if we consider the 
reparametrization 
$\Delta_\alpha=|\alpha_1-\alpha_2|$ and $\Sigma_\alpha=\alpha_1+\alpha_2$. As is the case for the 
Dirichlet distribution, the 
level of dependence is higher for large values of $\Sigma_\alpha$. Furthermore, $\lambda_u$ is 
monotonically decreasing in 
$\rho$.
Higher levels of extremal asymmetry, as measured by departures from the diagonal on the copula 
scale, are governed by both
$\Sigma_\alpha$ and $\Delta_\alpha$. The larger $\Sigma_\alpha$, the lower the asymmetry. Likewise, 
the larger $\Delta_\alpha$, 
the larger the asymmetry. 
Contrary to the case of extremal dependence, the behavior in $\rho$ is not monotone. 
For the negative scaled extremal Dirichlet model, asymmetry is 
maximal when $\rho \approx \alpha_1 \wedge \alpha_2$. When $\Sigma_\alpha$ is small, smaller values 
of $\rho$ induce larger 
asymmetry, but this is not the case for larger values of $\Sigma_\alpha$ where the asymmetry profile 
is convex with a global 
maximum attained for larger values of $\rho$.

\section{de Haan representation and simulation algorithms} \label{sec:6}

Random samples from the scaled extremal Dirichlet model can be drawn 
efficiently using the algorithms recently developed in \cite{Dombry:2016}. We first derive the 
so-called de Haan representation 
in 
Section \ref{sec:6.1} and adapt the algorithms from \cite{Dombry:2016} to the present setting in 
Section \ref{sec:6.2}.

\subsection{de Haan representation} \label{sec:6.1}
First, introduce the following family of univariate distributions, which we term the scaled Gamma 
family and denote by
$\mathsf{sGa}(a, b, 
c)$. It has three parameters $a, c > 0$ and $b \neq 0$ and a density given, for all $x >0$, by
\begin{equation}\label{eq:sGamma}
f(x; a, b, c) = \frac{|b|}{\Gamma(c)}a^{-bc}x^{bc-1} \exp \left\{ - \pfrac{x}{a}^b\right\}.
\end{equation}
Observe that when $Z \sim \mathsf{Ga}(c, 1)$ is a Gamma variable with shape parameter $c > 0$ and 
scaling parameter $1$, $Y 
\eqdis 
a Z^{1/b}$ is scaled Gamma $\mathsf{sGa}(a, b, c)$. Consequently, $\E{Y} = a 
\Gamma(c+1/b)/\Gamma(c) < \infty$ provided that $b < -1/c$. The scaled Gamma family includes several 
well-known distributions as 
special cases, notably the Gamma when $b=1$, the Weibull when $c=1$ and $b > 0$, the inverse Gamma 
when $b=-1$, and the 
Fr\'echet
when $c=1$ and $b < 0$. When $b > 0$, the scaled Gamma is the generalized Gamma distribution of 
\cite{Stacy:1962}, albeit in a 
different parametrization. 

Now consider the parameters $\bs{\alpha} =(\alpha_1, \dots, \alpha_d)$ with $\bs{\alpha} >\bs{0}_d$ 
and $\rho >- 
\min(\alpha_1, \ldots, \alpha_d)$, $\rho \neq 0$. Let $\bs{V}$ be a random vector with independent 
scaled 
Gamma margins $V_i \sim 
\mathsf{sGa}\{1/c(\alpha_i, \rho), 1/\rho, \alpha_i\}$, where for $\alpha >0$, $c(\alpha, \rho) = 
\Gamma(\alpha+\rho)/\Gamma(\alpha)$ 
as 
in Definition \ref{def:1}.  If $\bs{Z}$ is a random vector with independent Gamma margins $Z_i 
\sim 
\mathsf{Ga}(\alpha_i, 1)$ 
then 
for all $i=1, \dots, d$, $V_i \eqdis Z_i^{\rho}/c(\alpha_i, \rho)$. Furthermore, recall that 
$\|\bs{Z}\| \sim 
\mathsf{Ga}(\bar \alpha, 1)$ is independent of $\bs{Z}/\| \bs{Z}\|$, which has the same distribution 
as the Dirichlet vector 
$\bs{D}_{\bs{\alpha}}=(D_1, \dots, D_d)$.  One thus has, for all $\bs{x} \in \mathbb{R}_+^d$, 
\begin{align}\label{eq:deH}
\mathrm{E}\left\{\max_{1\le i \le d} (x_i V_i)\right\} = \mathrm{E}\left[\max_{1\le i \le d} \left\{ 
\frac{x_i 
Z_i^\rho}{c(\alpha_i, \rho)}\right\}\right] = \E{\|\bs{Z}\|^\rho} \mathrm{E}\left[\max_{1\le i \le 
d} \left\{ \frac{x_i 
D_i^\rho}{c(\alpha_i, \rho)}\right\}\right]=
\ell^{\mathrm{D}}(\bs{x};\rho, \bs{\alpha}), 
\end{align}
where $\ell^{\mathrm{D}}$ is as in \Cref{def:1}, given that $\E{\|\bs{Z}\|^\rho}=c(\bar \alpha, 
\rho)$.

When $\rho=1$, the positive scaled Dirichlet extremal model becomes the Coles--Tawn Dirichlet 
extremal model, $V_i \sim \mathsf{Ga}(\alpha_i, 1)$ and \Cref{eq:deH} reduces to the 
representation derived in \cite{Segers:2012b}. When $\bs{\alpha} = \bs{1}_d$, $\ell^{\mathrm{D}}$ 
becomes the stable tail 
dependence 
function  of the negative logistic model, $V_i$ is Weibull and \Cref{eq:deH} is the representation 
in Appendix~A.2.4 of 
\cite{Dombry:2016}. 
Similarly, when $\rho < 0$ and $\bs{\alpha} = \bs{1}_d$, the negative scaled Dirichlet extremal 
model becomes the logistic 
model, 
$V_i$ is Fr\'echet and \Cref{eq:deH} is the representation 
in Appendix~A.2.4 of \cite{Dombry:2016}. The requirement that $\rho > -\min(\alpha_1, \dots, 
\alpha_d)$ ensures that the 
expectation of $V_i$ is finite for all $i\in\{1, \ldots, d\}$.

\Cref{eq:deH} implies that the 
max-stable random vector $\bs{Y}$ with unit Fr\'echet margins and extreme-value copula with stable 
tail dependence function 
$\ell^{\mathrm{D}}(\cdot ; \rho, \bs{\alpha})$ admits the de Haan \citep{deHaan:1984} spectral 
representation
\begin{equation}\label{eq:deHpp}
\bs{Y} \eqdis \max_{k \in \mathbb{N}} \zeta_k\bs{V}_k, 
\end{equation}
where $\mathcal{Z}=\{\zeta_k\}_{k=1}^\infty$ is a Poisson point process on $(0, \infty)$ 
with intensity $\zeta^{-2}\d \zeta$ and $\bs{V}_k$ is an i.i.d. sequence of random vectors 
independent of $\mathcal{Z}$. 
Furthermore, the univariate margins of $\bs{V}_k$ are independent and such that $V_{kj} \sim 
\mathsf{sGa}\{1/c(\alpha_j, \rho), 1/\rho, \alpha_j\}$ for $j=1, \dots, d$ with 
$\E{\bs{V}_k}=\bs{1}_d$ for all $k\in \N$.

\subsection{Unconditional simulation} \label{sec:6.2}
The de Haan representation \eqref{eq:deHpp} offers, among other things, an easy route to 
unconditional simulation of max-stable 
random vectors that follow the scaled Dirichlet extremal model, as laid out in \cite{Dombry:2016} in 
the more general context of 
max-stable processes. To see how this work applies in the present setting, fix an arbitrary 
$j_0\in\{1, \dots, d\}$ and recall  
that the $j_0$th extremal function $\phi^+_{j_0}$ is given, almost surely, as $\zeta_k\bs{V}_k$ such 
that $Y_{j_0} = \zeta_k 
V_{kj_0}$.  
From \cref{eq:deHpp} and Proposition~1 in \cite{Dombry:2016} it then directly 
follows that 
$\phi^+_{j_0} / Y_{j_0} \eqdis ( W_{j_01}/W_{j_0 j_0}, \ldots, W_{j_0 d}/W_{j_0 j_0})$, 
where $\bs{W}_{j_0} = (W_{j_0 1}, \dots, W_{j_0 d})$ is a random vector with density given, for all 
$\bs{x} \in 
\mathbb{R}_+^d$, 
by 
\begin{align*}
 \frac{|1/\rho|}{\Gamma(\alpha_{j_0})}c(\alpha_{j_0}, 
\rho)^{\alpha_{j_0}/\rho}x_{j_0}^{\alpha_{j_0}/\rho} \exp \left[ - 
\{c(\alpha_{j_0}, \rho) x_{j_0}\}^{1/\rho}\right] \times\prod_{j=1, j\neq j_0 }^d 
\frac{|1/\rho|}{\Gamma(\alpha_j)}c(\alpha_j, 
\rho)^{\alpha_j/\rho}x_j^{\alpha_j/\rho-1} \exp \left[ - \{c(\alpha_j, \rho) x_j\}^{1/\rho}\right]. 
\end{align*}
This means that the components of $\bs{W}_{j_0}$ are independent and such that $W_{j_0 j} \sim 
\mathsf{sGa}\{1/c(\alpha_j, \rho), 1/\rho, \alpha_j\}$ when $j \neq j_0$ and $W_{j_0 j_0} \sim 
\mathsf{sGa}\{1/c(\alpha_{j_0}, \rho), 1/\rho, \alpha_{j_0}+\rho\}$. In other words, $W_{j_0 j_0} 
\sim 
Z_{j_0}^\rho/c(\alpha_{j_0} , 
\rho)$ where $Z_{j_0} \sim \mathsf{Ga}(\alpha_{j_0}+\rho, 1)$, while for all $j \neq j_0$, $W_{j_0 
j}\eqdis 
Z_j^\rho/c(\alpha_j, 
\rho) $ where $Z_j \sim 
\mathsf{Ga}(\alpha_j, 1)$.

The exact distribution of $\phi^+_{j_0} / Y_{j_0}$ given above now allows for an easy adaptation of 
the algorithms 
in \cite{Dombry:2016}. To draw an observation from the extreme-value copula with the scaled 
Dirichlet stable tail dependence 
function $\ell^{\mathrm{D}}$ with parameters $\bs{\alpha} >\bs{0}_d$ and $\rho > -\min(\alpha_1, 
\dots, \alpha_d)$, $\rho \neq 
0$, one 
can follow Algorithms~1 and~2 below. The first procedure corresponds to Algorithm~1 in  
\cite{Dombry:2016} and relies on  
\cite{Schlather:2002}; the second is an adaptation of Algorithm~2 in \cite{Dombry:2016}.

\begin{algorithm}
\caption{Exact simulations from the extreme-value copula based on spectral densities.}
\label{algo:1}
\small
\begin{algorithmic}[1]
\State Simulate $E \sim \mathsf{Exp}(1)$ .
\State Set $\bs{Y} = \bs{0}$.
\While{$1/E > \min(Y_1, \dots, Y_d)$}
\State Simulate $J$ from the uniform distribution on $\{1, \dots, d\}$.
\State Simulate independent $Z_{j} \sim \mathsf{Ga}(\alpha_j, 1)$ for $j\in \{1, \ldots, 
d\}\setminus J$ and $Z_{J} \sim 
\mathsf{Ga}(\alpha_J + \rho, 1)$.
\State Set $W_j \gets Z_j^{\rho}/c(\alpha_j, \rho)$, $j=1, \dots, d$.
\State Set $\bs{S} \gets \bs{W}/\|\bs{W}\|$.
\State Update $\bs{Y} \gets \max\{\bs{Y}, d\bs{S}/E\}$.
\State Simulate $E^* \sim \mathsf{Exp}(1)$ and update $E \gets E+E^*$.
\EndWhile
\State \textbf{return} $\bs{U} = \exp(-1/\bs{Y})$.
\end{algorithmic}
\end{algorithm}
\begin{algorithm}
\caption{Exact simulations based on sequential sampling of the extremal functions.}
\label{algo:2}
\small
\begin{algorithmic}[1]
\State Simulate $Z_{1} \sim \mathsf{Ga}(\alpha_1 + \rho, 1)$ and $Z_{j} \sim 
\mathsf{Ga}(\alpha_j, 1)$, $j =2, \dots, d$.
\State Compute $\bs{W}$ where $W_{j} \gets Z_j^{\rho}/c(\alpha_j, \rho)$, $j=1, \dots, d$. 
\State Simulate $E_1 \sim \mathsf{Exp}(1)$.
\State Set $\bs{Y} \gets \bs{W}/(W_1E_1)$.
\For{$k=2, \dots, d$}
\State Simulate $E_k \sim \mathsf{Exp}(1)$.
\While{$ 1/E_k > Y_k$}
\State Simulate independent $Z_{k} \sim \mathsf{Ga}(\alpha_k + \rho, 1)$ and $Z_{j} \sim 
\mathsf{Ga}(\alpha_j, 1)$, $j =1, \ldots, d, j\neq k$ .
\State Set $\bs{W}=(W_1, \ldots, W_d)$ where $W_{j} \gets Z_j^{\rho}/c(\alpha_j, \rho)$, $j=1, 
\dots, d$. 
\If{$W_{i}/(W_kE_k) < Y_i$ for all $i=1, \dots, k-1$}
\State Update $\bs{Y} \gets \max\{\bs{Y}, \bs{W} /(W_kE_k)\}$.
\EndIf
\State Simulate $E^* \sim \mathsf{Exp}(1)$ and update $E_k \gets E_k + E^*$.
\EndWhile
\EndFor
\State \textbf{return} $\bs{U} = \exp(-1/\bs{Y})$.
\end{algorithmic}
\end{algorithm}

Note that $\bs{S}$ obtained in Step~{\footnotesize{7}} of Algorithm~1 has the angular distribution 
$\sigma_d$ of 
$\ell^{\mathrm{D}}$; 
see Theorem~1 in \cite{Dombry:2016}. Similar algorithms for drawing samples from the angular 
distribution of the extremal logistic
and Dirichlet models were obtained in \cite{Boldi:2009}. Algorithm~2 requires a lower number of 
simulations and is 
more efficient on average, cf. \cite{Dombry:2016}.  Both algorithms are easily implemented using the 
function \texttt{rmev} in 
the 
\texttt{mev} package within the \textsf{R} Project for Statistical Computing \cite{Rlang}, which 
returns samples of max-stable scaled extremal Dirichlet vectors with unit Fr\'echet margins, i.e., 
$\bs{Y}$ in \Cref{algo:1,algo:2}.

\section{Estimation} \label{sec:7}

The scaled extremal Dirichlet model can be used to model dependence between extreme events. To this 
end, 
several schemes can be envisaged. For example, one can consider the block-maxima approach, given 
that max-stable distributions 
are  the most natural for such data. Another option is peaks-over-threshold models.
Yet another alternative, used in \cite{Engelke:2015} for the Brown--Resnick model, is to approximate 
the conditional 
distribution of a 
random vector with unit Fr\'echet margins given that the $j_0$th component exceeds a large threshold 
by the 
distribution of 
$\phi^+_{j_0}/Y_{j_0}$ discussed in Section~\ref{sec:6.2}.

Here, we focus on the multivariate tail model of \cite{Ledford:1996}; see also Section~16.4 in 
\cite{McNeil/Frey/Embrechts:2015}. 
To this end, let $\bs{X}_1, \dots, \bs{X}_n$ be a random sample from some unknown multivariate 
distribution $H$ with continuous 
univariate margins which is assumed to be in the maximum domain of attraction of a multivariate 
extreme-value 
distribution $H_0$. To model 
the tail of $H$, its margins $F_j$, $j=1, \dots, d$ can first be approximated using the univariate 
peaks-over-threshold 
method. For 
all $x$ above some high threshold $u_j$, one then has
\begin{align}\label{eq:margpareto}
 F_j(x) \approx  \tilde F_j(x; \eta_j, \xi_j)= 1- 
\nu_j\left(1+\xi_j\frac{(x-u_j)}{\eta_j}\right)_{+}^{-1/\xi_j}, 
\end{align}
where $\nu_j = 1-F_j(u_j)$, and $\eta_j > 0$ and $\xi_j $ are the parameters of the generalized 
Pareto distribution. Furthermore, 
for $\bs{w}$ sufficiently close to $\bs{1}_d$, the copula of $H$ can be approximated by the 
extreme-value copula $C_0$ of $H_0$, 
so that, for $\bs{x} \ge \bs{u}$, $H(\bs{x}) \approx \tilde H(\bs{x}) =C_0\{\tilde F_1(x_1), \dots, 
\tilde F_d (x_d)\}$. The 
parameters of this multivariate tail model, i.e., the parameters $\bs{\theta}$ of the stable tail 
dependence function $\ell_0$ 
of $C_0$ as well as the marginal parameters $\bs{\nu}$, $\bs{\eta}$ and $\bs{\xi}$ can be estimated 
using likelihood methods; 
this  allows, e.g., for Bayesian inference, generalized additive modeling of the parameters and 
model selection based on 
likelihood-ratio tests. For a comprehensive review of likelihood inference methods for 
extremes, see, e.g.,  \cite{Huser:2016}.

The multivariate tail model can be fitted in low-dimensions using the censored likelihood $L(\bs{X}; 
\bs{\nu}, 
\bs{\eta}, \bs{\xi}, \bs{\theta}) = \prod_{i=1}^n  L_i (\bs{X}_i; \bs{\nu}, \bs{\eta}, \bs{\xi}, 
\bs{\theta})$, where for 
$i=1, \dots, n$, 
\begin{align}\label{eq:clik}
L_i (\bs{X}_i; \bs{u}, \bs{\nu}, \bs{\eta}, \bs{\xi}, \bs{\theta})  = \left.\frac{\partial^{m_i} 
\tilde 
H(y_1, \dots, y_d)}{\partial y_{j_1} \dotsm \partial y_{j_{m_i}}}\right|_{\bs{y}=\max(\bs{X}_i, 
\bs{u})} 
= \left.\frac{\partial^{m_i} \exp\{-\ell_0 (1/\bs{y})\}}{\partial y_{j_1} \dotsm \partial 
y_{j_{m_i}}}\right|_{\bs{y}=t\{\max(\bs{X}_i, \bs{u})\}} 
\prod_{k=1}^{m_i} J_{j_k}(X_{ij_k})
\end{align}
In this expression, the indices $j_1, \dots, j_{m_i}$ are those of the components of $\bs{X}_i$ 
exceeding the thresholds 
$\bs{u}$ and  for $\bs{x} \ge \bs{u}$, $t(\bs{x})  = (t_1(x_1), \dots, t_d(x_d))$, where for $j=1, 
\dots, d$, 
\begin{align}\label{eq:tJ}
t_j(x_j) = -\frac{1}{\log \{\tilde{F}_{j}(x_j;\eta_j, \xi_j)\}}, 
\qquad  
J_j(x_j)=\frac{\nu_{j}}{\eta_{j}} \left(1+\xi_j\frac{(x_j-u_j)}{\eta_j}\right)^{-1/\xi_j-1} 
\!\!\!\frac{1}{[\log 
\{\tilde{F}_{j}(x_j;\eta_j, \xi_j)\}]^2 \tilde{F}_{j}(x_j;\eta_j, \xi_j)}.
\end{align}
The censored likelihood $L(\bs{X}; \bs{\nu}, 
\bs{\eta}, \bs{\xi}, \bs{\theta})$ can be maximized either over all parameters at once, or the 
marginal parameters $\bs{\nu}$, 
$\bs{\eta}$ and $\bs{\xi}$ can be estimated from each univariate margin separately, so that only the 
estimate 
of $\bs{\theta}$ is obtained through maximizing $L$. When $d$ is large, one can also maximize the 
likelihood in 
\cite{Smith:1997} that uses the tail approximation $\bar{H}(\bs{x}) \approx 1-\ell(1/\bs{x})$. In 
either case, 
$\ell_0$ and the higher-order partial derivatives of $\ell_0(1/\bs{x})$ need to be computed.

When $\ell_0$ is the scaled 
extremal Dirichlet stable tail dependence function $\ell^{\mathrm{D}}(\cdot ; \rho, \bs{\alpha})$ 
given in \Cref{def:1} with 
parameters $\bs{\alpha} > 0$ and $\rho > -\min(\alpha_1, \dots, \alpha_d)$, $\rho \neq 0$, its 
expression is not explicit. 
However, 
$\ell^{\mathrm{D}}$ can be calculated numerically using adaptive numerical cubature algorithms for 
integrals of functions defined 
on 
the simplex, as implemented in, e.g., the \textsf{R} package 
\href{https://cran.r-project.org/web/packages/SimplicialCubature/index.html}{\texttt{SimplicialCubature}}. 
Given the 
representation in \cref{eq:ellD}, $\ell^{\mathrm{D}}$ is also easily approximated using Monte Carlo 
methods. Instead of employing 
\cref{eq:ellD} directly 
and sampling from the Dirichlet vector $\bs{D}_{\bs{\alpha}}$, one can use the more efficient 
importance sampling 
estimator
\begin{align*}
\widehat{\ell^{\mathrm{D}}}(1/\bs{u}, \rho, \bs{\alpha})= \frac{1}{B} \sum_{i=1}^B \frac{\max_{1 \le 
j\le d} \left[
\{c(\alpha_j, \rho)u_j\}^{-1}D_{ij}^\rho\right]}{\frac{1}{d}\sum_{j=1}^d {c(\alpha_j, 
\rho)^{-1}D_{ij}^\rho}}, 
\end{align*}
where $\bs{D}_{i} \sim d^{-1} \sum_{j=1}^d \mathsf{Dir}(\bs{\alpha}+\I{j}\rho\bs{1}_d)$ is sampled 
from a
Dirichlet mixture. 

The partial  derivatives of $\ell^{\mathrm{D}}$ can be calculated using the following result, 
shown in \ref{app:E}.
\begin{proposition} \label{prop:6}
Let $\ell^{\mathrm{D}}$ be the scaled extremal Dirichlet stable tail dependence function with 
parameters $\bs{\alpha} 
>\bs{0}_d$ 
and $-\min(\alpha_1, \dots, \alpha_d) < \rho < \infty$, $\rho \neq 0$. Then, for any $\bs{x} \in 
\mathbb{R}_+^d$, 
\begin{align}
\frac{\partial^d \ell^{\mathrm{D}}(1/\bs{x})}{\partial x_1 \cdots \partial 
x_d}&= - d h^{\mathrm{D}}(\bs{x};\rho, \bs{\alpha})= 
-\frac{\Gamma(\bar{\alpha}+\rho)}{|\rho|^{d-1}\prod_{i=1}^d\Gamma(\alpha_i)}
\left[\sum_{j=1}^d \left\{c(\alpha_j, 
\rho)x_j\right\}^{1/\rho}\right]^{-\rho-\bar{\alpha}}\prod_{i=1}^{d} 
\{c(\alpha_i, \rho)\}^{\alpha_i/\rho}x_i^{\alpha_i/\rho-1}, 
\label{eq:densitypD}
\end{align}
where $h^{\mathrm{D}}$ is as given in Proposition \ref{prop:5}.
Furthermore, for all $k=1, \dots, d-1$ and $\bs{x} \in \mathbb{R}_+^d$, 
\begin{align*}
\frac{\partial^k \ell^{\mathrm{D}}(1/\bs{x})}{\partial x_1 \cdots \partial x_k} = - \int_0^\infty 
t^{k} \prod_{i=1}^k f
\left(x_i t; \frac{1}{c(\alpha_i, \rho)}, \frac{1}{\rho}, \alpha_i\right)\prod_{i=k+1}^d 
F\left(x_i t; \frac{1}{c(\alpha_i, \rho)}, \frac{1}{\rho}, \alpha_i\right) \d t, 
\end{align*}
where $f(; a, b, c)$ and $F(;a, b, c)$ denote, respectively the density and distribution function of 
the scaled Gamma 
distribution 
with parameters $a, c > 0$ and $b \neq 0$ given in \cref{eq:sGamma}. Furthermore, if $\gamma(c, x) = 
\int_0^x t^{c-1} e^{-t} dt$ 
denotes the lower incomplete gamma function, then for $x >0$, $F(x;a, b, c) = \gamma\{ c, (x/a)^b\} 
/\Gamma(c)$ when $b > 0$ 
while 
$F(x;a, b, c) = 1-\gamma\{ c, (x/a)^b\} /\Gamma(c)$ when $b < 0$.
\end{proposition}
Other estimating equations could be used to circumvent the calculation of 
$\ell^{\mathrm{D}}(1/\bs{x})$ and its partial 
derivatives. 
An interesting alternative to likelihoods in the context of proper 
scoring functions is proposed in \cite{deFondeville:2016}. Specifically, the authors 
advocate the use of the gradient score, adapted by them for the peaks-over-threshold framework, 
\begin{align*}
  \delta_w(\bs{x})=\sum_{i=1}^d \left(2w_i(\bs{x})\frac{\partial w_i(\bs{x})}{\partial 
x_i}\frac{\partial 
\log 
h(\bs{x})}{\partial x_i} + w_i^2(\bs{x}) \left[ \frac{\partial^2\log h(\bs{x})}{\partial x_i^2} + 
\frac{1}{2}\left\{\frac{\partial \log h(\bs{x})}{\partial x_i}\right\}^2\right]\right)
\end{align*}
for a differentiable weighting function $w(\bs{x})$, unit Fr\'echet observations $\bs{x}$ and 
density $h(\bs{x})$ that 
would 
correspond in the setting of the scaled extremal Dirichlet to $d h^{\mathrm{D}}(\bs{x};\rho, 
\bs{\alpha})$. Explicit expressions 
for 
the derivatives of $\log d h^{\mathrm{D}}$ may be found in \ref{app:E}. The 
parameter estimates are obtained as the solution to 
$\argmax_{\bs{\theta} \in \Theta} \sum_{i=1}^n\delta_w(\bs{x}_i)\I{\mathscr{R}(\bs{x}_i/\bs{u})>1}$, 
where 
$\bs{\theta}=(\rho, \bs{\alpha})$ is the vector of parameters of the model and $\mathscr{R}$ is a 
differentiable risk functional, 
usually the $\ell_p$ norm for some 
$p \in \N$.
Although the gradient score is not asymptotically most efficient, weighting functions can be 
designed to reproduce approximate 
censoring, lending the method robustness and tractability.

\section{Data illustration} \label{sec:8}
\begin{figure}[t!]
\centering 
\includegraphics[width=0.9\textwidth]{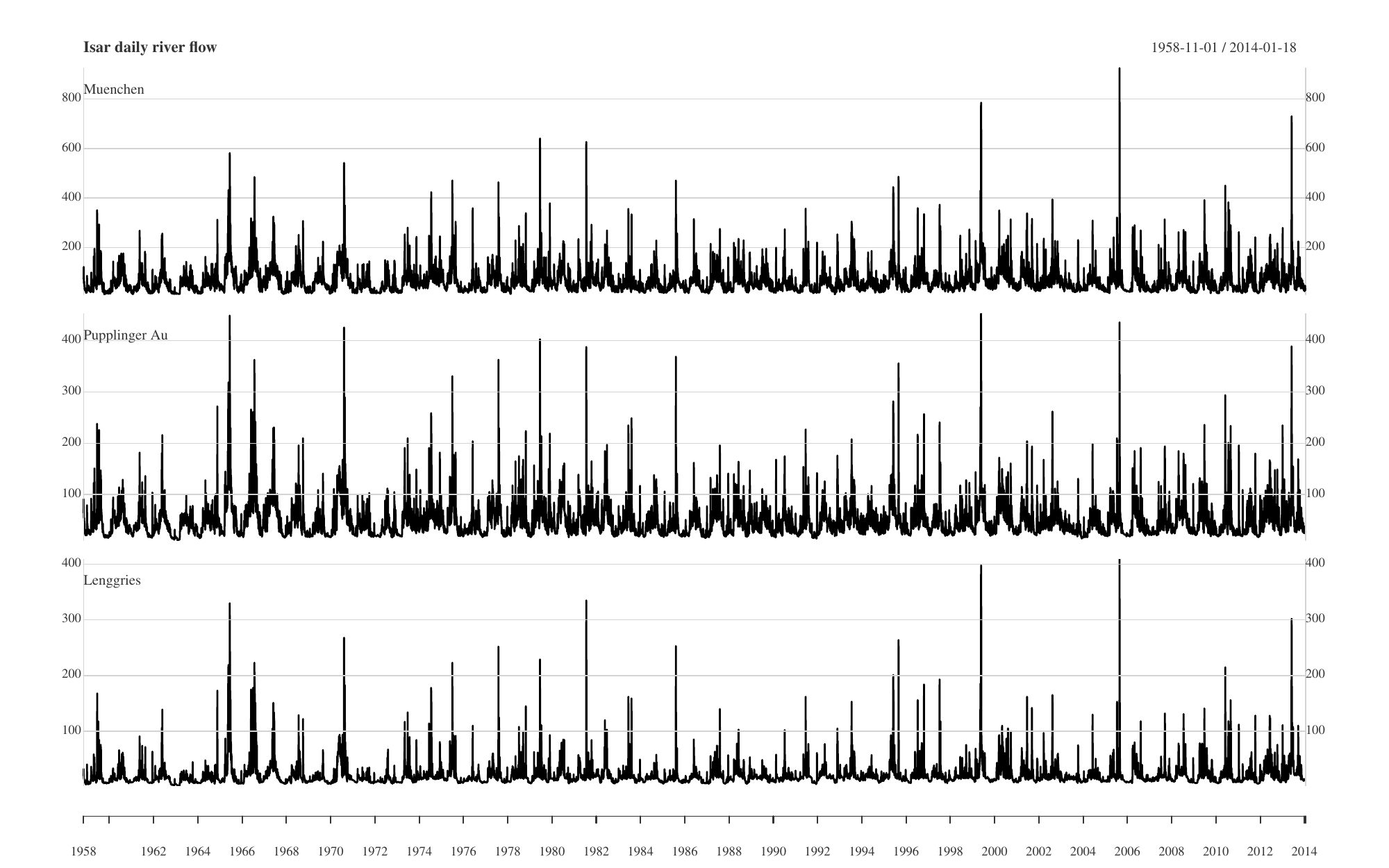}
\caption{Daily river flow of the Isar river at three sites}
\label{fig:data}
\end{figure}

In this section, we illustrate the use of the scaled extremal Dirichlet model on a trivariate sample 
of daily river flow data of 
the river Isar in southern  Germany; this dataset is a subset of the one analyzed in 
\cite{Asadi:2015}. All the code can be 
downloaded from \href{https://github.com/lbelzile/ealc}{\texttt{https://github.com/lbelzile/ealc}}. For this analysis, we 
selected data measured 
at Lenggries (upstream), 
Pupplinger Au (in the middle) and Munich (downstream). To ensure stationarity of the series and given 
that the most extreme events 
occur during the summer, we restricted our attention to the measurements for the months of June, 
July and August. Since
the sites are measuring the flow of the same river, dependence at extreme levels is likely to be 
present, as is indeed apparent 
from \Cref{fig:data}.  
Directionality of the river may further lead to asymmetry in the asymptotic dependence structure, 
suggesting that the scaled  
extremal Dirichlet model may be well suited for these data. Furthermore, given that other well-known 
models like the extremal 
Dirichlet, logistic and negative logistic are nested within this family, their adequacy can be 
assessed through 
likelihood ratio tests.

To remove dependence between extremes over time, we decluster each series and retain only the 
cluster maxima based on three-day 
runs. 
Rounding of the measurements has no impact on parameter estimates and is henceforth neglected. 
The multivariate tail model outlined in Section~\ref{sec:7} is next fitted to the cluster maxima. 
The thresholds 
$\bs{u}=(u_1, u_2, u_3)$  were selected to be the 92\%  quantiles using the parameter stability plot 
of \cite{Wadsworth:2016} 
(not 
shown here). Next, set $\bs{\theta}=(\bs{\eta}, \bs{\xi}, \bs{\alpha}, \rho)$, where $\bs{\eta}$ and 
$\bs{\xi}$ are the 
marginal 
parameters of the generalized Pareto distribution in \cref{eq:margpareto} and $\rho$ and 
$\bs{\alpha}$ are the parameters of the 
scaled Dirichlet model. To estimate $\bs{\theta}$, the trivariate censored likelihood 
\eqref{eq:clik} could be used. To avoid 
numerical integration and because of the relative robustness to misspecification, we employed the 
pairwise composite 
log-likelihood $l_C$ of \cite{Ledford:1996} instead;  the loss of efficiency in this trivariate 
example is likely small. 
Specifically, we maximized
\begin{align*}
  l_C(\bs{\theta})=\sum_{i=1}^n \sum_{j=1}^{d-1} \sum_{k=j+1}^d 
  \left[ \log g\{t_j(x_{ij}), t_k(x_{ik}); \bs{\theta}, t_j({u}_j), t_k(u_k)\} + 
\I{x_{ij}>u_{j}}\log J_j(x_{ij}) + \I{x_{ik}>u_{k}}\log J_k(x_{ik})\right], 
\end{align*}
where
\begin{align*}
  g(y_{j}, y_{k}; \bs{\theta}, u_j, u_k) = 
\begin{cases}
\exp\{-\ell(1/u_j, 1/u_k)\}, &  y_j \leq u_j, y_k \leq u_k\\
-{\partial \ell(1/y_j, 1/u_k)/\partial y_j}\exp\{-\ell(1/y_j, 1/u_k)\}, &  y_j > u_j, y_k \leq u_k\\
-{\partial \ell(1/u_j, 1/y_k)}/{\partial y_k} \exp\{-\ell(1/u_j, 1/y_k)\}, &  y_j \leq u_j, y_k > 
u_k\\
\left[\left\{\partial \ell(1/y_j, 1/y_k)/{\partial y_j}\right\} \left\{{\partial \ell(1/y_j, 
1/y_k)}/{\partial 
y_k}\right\}-dh^{\mathrm{D}}(y_j, y_k)\right]\exp\{-\ell(1/y_j, 1/y_k)\}, &  y_j > u_j, y_k > u_k
\end{cases}
\end{align*}
where $\ell=\ell^{\mathrm{D}}$ and for all $j=1, \dots, d$, $t_j$ and $J_j$ are as in \Cref{eq:tJ}.

Uncertainty assessment can be done in the same way as for general estimating equations. 
Specifically, let $g(\bs{\theta})$ 
denote an unbiased 
estimating function and define the variability matrix $\mathbf{J}$, 
the sensitivity matrix $\mathbf{H}$ and the Godambe information matrix $\mathbf{G}$ as
\begin{align}
\mathbf{J}=\E{\frac{\partial g(\bs{\theta})}{\partial \bs{\theta}}\frac{\partial 
g(\bs{\theta})}{\partial 
\bs{\theta}}^\top}, \qquad 
\mathbf{H}=-\E{\frac{\partial^2 g(\bs{\theta})}{\partial \bs{\theta}\partial \bs{\theta}^\top}}, 
\qquad 
\mathbf{G}=\mathbf{H}\mathbf{J}^{-1}\mathbf{H}. \label{eq:godambe}      
\end{align}
The maximum composite likelihood estimator is strongly consistent and asymptotically normal, 
centered at the true parameter $\bs{\theta}$ with covariance matrix given by the inverse Godambe 
matrix $\mathbf{G}^{-1}$.

\begin{figure}[t!]
\centering 
\includegraphics[width=0.8\textwidth]{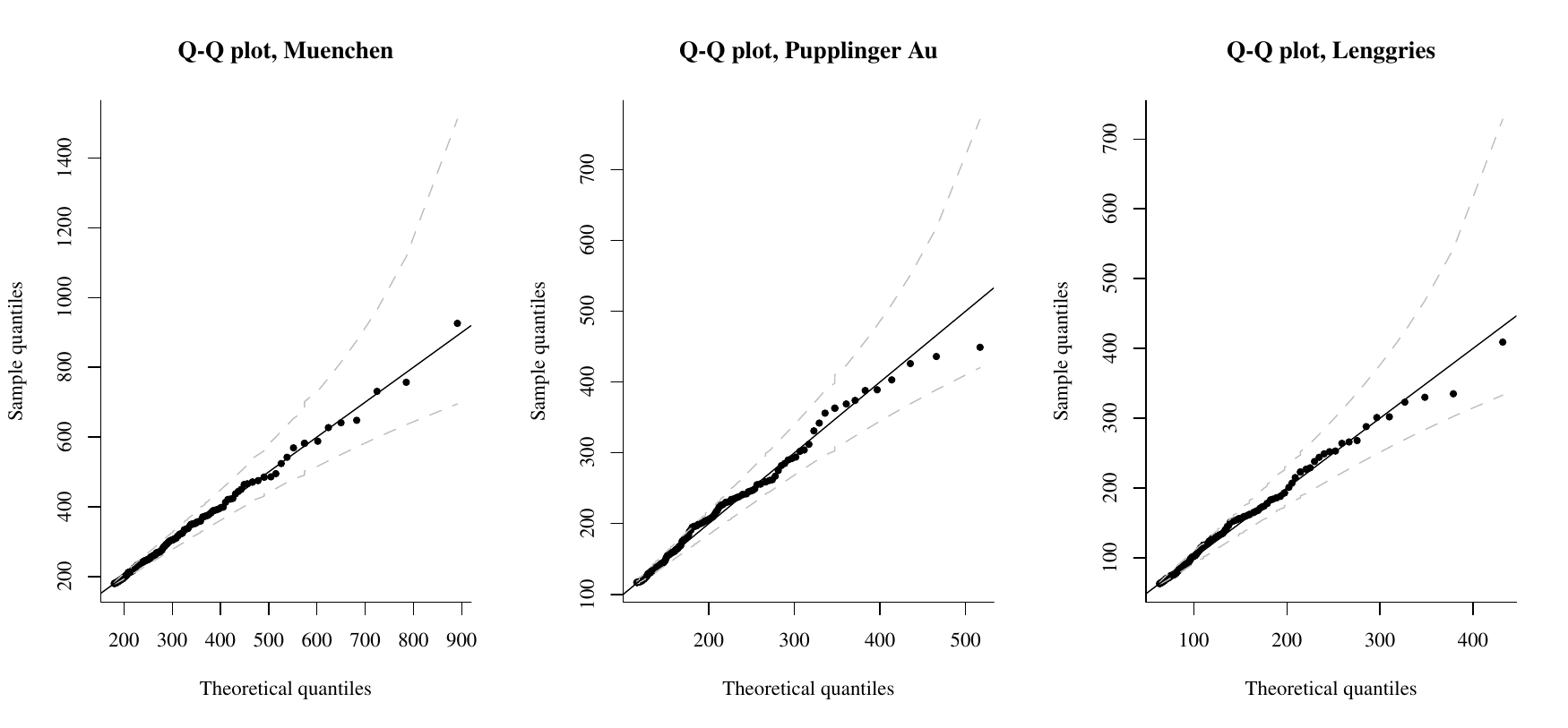}
\caption{Marginal Q-Q plots for the three sites based on pairwise composite likelihood estimates for 
the scale and shape 
parameters obtained from the scaled Dirichlet model, retaining marginal exceedances of the 92\% 
quantiles. The pointwise 
confidence intervals were obtained from the 
transformed Beta quantiles of the order statistics.} 
\label{fig:qqplot}
\end{figure}

Using the pairwise composite log-likelihood $l_C$, we fitted the scaled extremal Dirichlet model as 
well as the logistic and 
negative logistic models that correspond to the negative and positive scaled extremal Dirichlet 
models, respectively, and the 
parameter restriction $\bs{\alpha} = \bs{1}_d$.  The estimates of the marginal generalized Pareto 
parameters $\bs{\eta}$ and 
$\bs{\xi}$ are given in \Cref{tab:margpar}. 	
As the estimates were obtained by maximizing $\ell_C$, their values depend on the fitted model; the 
line labeled ``Marginal" 
corresponds to fitting the generalized Pareto distribution to threshold exceedances of each one of 
the three series separately. 
The marginal Q-Q plots displayed in \Cref{fig:qqplot} indicate a 
good fit of the model as well.

\begin{table}[ht]
\centering
\begin{tabular}{rllllll}
  \toprule
 & $\eta_1$ & $\eta_2$ & $\eta_3$ & $\xi_1$ & $\xi_2$ & $\xi_3$ \\ 
  \midrule
Scaled Dirichlet & 123.2 (7.5) & 84.4 (5) & 68.1 (4.2) & 0.05 (0.04) & -0.03 (0.04) & 0.02 (0.04) 
\\ 
  Neg. logistic & 117.1 (6.8) & 86.2 (5.1) & 70 (4.3) & 0.08 (0.04) & -0.05 (0.04) & 0 (0.04) \\ 
  Logistic & 117.3 (6.8) & 86.6 (5.1) & 70.4 (4.3) & 0.08 (0.04) & -0.05 (0.04) & 0 (0.04) \\ 
  ext. Dirichlet & 114.4 (6.8) & 84.3 (4.9) & 68.2 (4.1) & 0.12 (0.04) & -0.02 (0.04) & 0.04 (0.04) 
\\ 
  Marginal & 129.1 (14.5) & 95.1 (10.6) & 76 (8.7) & -0.01 (0.08) & -0.15 (0.08) & -0.08 (0.08) \\ 
   \bottomrule
\end{tabular}
\caption{Generalized Pareto parameter estimates and standard errors (in parenthesis) for the 
trivariate river example for four 
different models. } 
\label{tab:margpar}
\end{table}

\begin{table}[ht]
\centering
\begin{tabular}{rllll}
  \toprule
 & $\alpha_1$ & $\alpha_2$ & $\alpha_3$ & $\rho$ \\ 
  \midrule
Scaled Dirichlet & 0.76 (0.3) & 1.65 (0.82) & 2.03 (1.15) & $-$0.32 (0.1) \\ 
  Neg. logistic & 1 & 1 & 1 & 0.36 (0.02) \\ 
  Logistic & 1 & 1 & 1 & 0.28 (0.01) \\ 
  ext. Dirichlet & 3.34 (0.52) & 10.2 (2.84) & 12.78 (3.93) & 1 \\ 
  Gradient score & 1 & 2.72 & 2.66 & $-$0.39 \\ 
   \bottomrule
\end{tabular}
\caption{Dependence parameters estimates and standard errors (in parenthesis) for the trivariate 
river example.} 
\label{tab:deppar}
\end{table}

The estimates of the dependence parameters $\bs{\alpha}$ and $\rho$ are given in \Cref{tab:deppar}. 
The last line displays the maximum gradient score estimates were obtained from the raw data, i.e., 
ignoring the clustering,  
after 
transforming the observations to the standard Fr\'echet scale using the probability integral 
transform. We 
 retained only the 10\% largest values based on the $\ell_p$ norm with $p=20$; this risk functional 
is essentially a 
differentiable approximation of $\ell_\infty$.  We selected the weight function $w(\bs{x}, 
u)=\bs{x}[1 - 
\exp\{-(\|\bs{x}\|_p/ u - 1)\}]$ based on \cite{deFondeville:2016} to reproduce approximate 
censoring. The estimates are 
similar to the composite maximum likelihood estimators, though not efficient. 

\begin{figure}[t!]
\centering 
\includegraphics[width=0.9\textwidth]{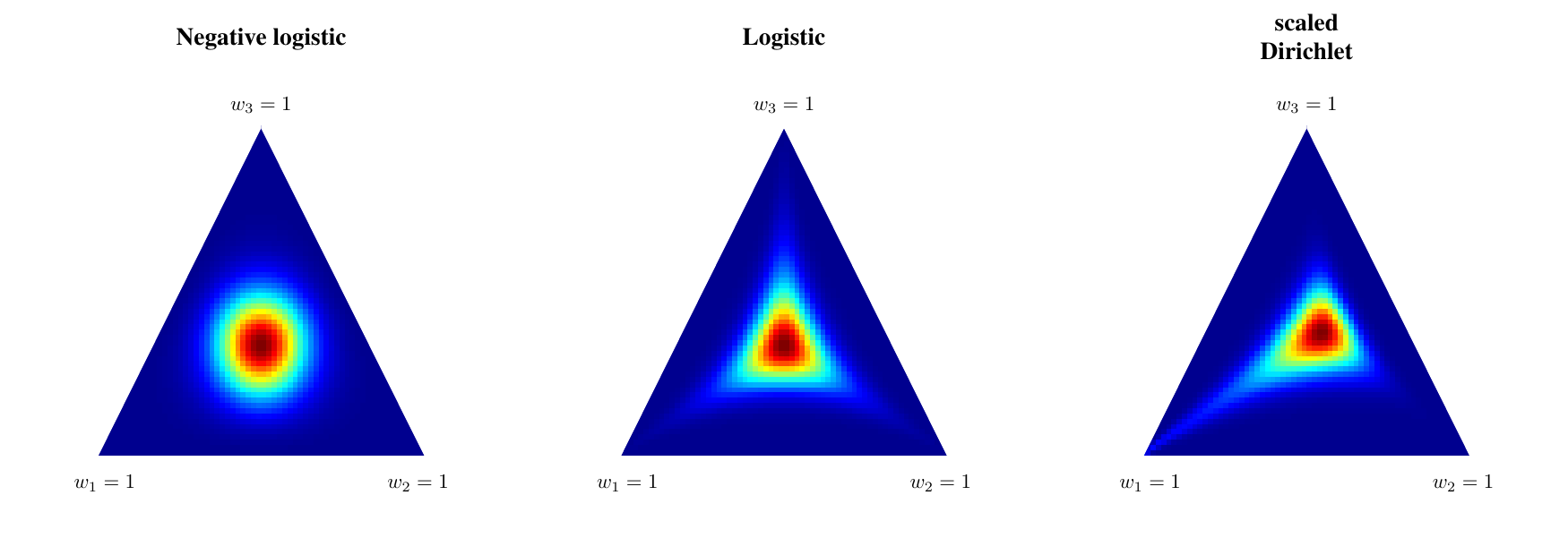}
\caption{Angular density plots for the three models, the negative logistic (left), logistic (middle) 
and scaled 
Dirichlet (right). The colours correspond to log density values and range from red (high density) 
to 
blue (low density). } 
\label{fig:angdens}
\end{figure}

The angular densities of the fitted logistic, negative logistic and scaled extremal Dirichlet models 
are displayed in 
\Cref{fig:angdens}. The right panel of this figure shows asymmetry caused by a few extreme events 
that only 
happened downstream. Whether this asymmetry is significant can be assessed through composite 
likelihood ratio tests; recall that 
the logistic model, the negative logistic model and the extremal Dirichlet model of \cite{Coles:1991} are all 
nested 
within the scaled extremal Dirichlet model. To this end, consider a partition of 
$\bs{\theta}=(\bs{\psi}, \bs{\lambda})$ into a 
$q$ dimensional parameter of interest $\bs{\psi}$ and a 
$3d+1-q$ dimensional nuisance parameter $\bs{\lambda}$, and the corresponding partitions of the 
matrices $\mathbf{H}$, 
$\mathbf{J}$ and $\mathbf{G}$.
Let $\widehat{\bs{\theta}}_C=(\widehat{\bs{\psi}}_C, \widehat{\bs{\lambda}}_C)$ denote the maximum 
composite likelihood 
parameter 
estimates and  $\widehat{\bs{\theta}}_0=(\bs{\psi}_{0}, \widehat{\bs{\lambda}}_{0})$ the restricted 
parameter estimates under 
the 
null hypothesis that the simpler model is adequate. The asymptotic distribution of the composite 
likelihood ratio test statistic 
$2\{\log l_C(\widehat{\bs{\theta}_C})-\log l_C(\widehat{\bs{\theta}_0})\}$ is
equal to $\sum_{i=1}^q c_iZ_i$
where $Z_i$ are independent $\chi^2_1$ variables and $c_i$ are the eigenvalues of the $q \times q$ 
matrix 
$(\mathbf{H}_{\bs{\psi\psi}}-\mathbf{H}_{\bs{\psi\lambda}}\mathbf{H}^{-1}_{\bs{\lambda\lambda}}
\mathbf{H}_{\bs{\lambda\psi}}
)\mathbf{G}^{-1}_{\bs{\psi\psi}}$; see \cite{Kent:1982}. 
We estimated the inverse Godambe information matrix, $\smash{\mathbf{G}^{-1}}$, by the empirical 
covariance of ${B}$ 
nonparametric bootstrap replicates.
The sensitivity matrix $\mathbf{H}$ was obtained from the Hessian matrix at the
maximum composite likelihood estimate and the variability matrix $\mathbf{J}$ from 
\cref{eq:godambe}. 
Since the Coles--Tawn extremal Dirichlet, negative logistic and logistic models are nested within 
the scaled Dirichlet family, 
we test for a restriction to these simpler models; the respective approximate $P$-values were 0.003, 
0.74 and 0.78. These values 
suggest that while the Coles--Tawn extremal Dirichlet model is clearly not suitable, there is not 
sufficient evidence 
to discard the logistic and negative logistic models. 
The effects of possible model misspecification are also visible for the Coles--Tawn extremal 
Dirichlet model, as the parameter 
values of $\alpha_1, \alpha_2$ and $\alpha_3$ are very large (viz. \Cref{tab:deppar}) and this 
induces negative
bias in the shape parameter estimates, as can be seen from \Cref{tab:margpar}.

\section{Discussion}\label{sec:9}

In this article, we have identified extremal attractors of copulas and survival copulas of Liouville 
random vectors $ R 
\bs{D}_{\bs{\alpha}}$, where $\bs{D}_{\bs{\alpha}}$ has a Dirichlet distribution on the unit simplex 
with parameters 
$\bs{\alpha}$, and $R$ is a strictly positive random variable independent of $\bs{D}_{\bs{\alpha}}$. 
The limiting stable tail 
dependence functions can be embedded in a single family, which can capture asymmetry and provides a 
valid model in dimension $d$. 
The 
latter is novel and termed here the scaled extremal Dirichlet; it 
includes the well-known logistic, negative logistic as well as the Coles--Tawn extremal Dirichlet 
models as special cases. In 
particular, therefore, this paper is first to provide an example of a random vector attracted to the 
Coles--Tawn extremal 
Dirichlet model, which was derived by enforcing moment constraints on a simplex distribution rather 
than as the limiting 
distribution of a random vector. 

A scaled extremal Dirichlet stable tail dependence function $\ell^{\mathrm{D}}$ has $d+1$ 
parameters, $\rho$ and $\bs{\alpha}$. 
The 
parameter vector $\bs{\alpha}$ is inherited from $\bs{D}_{\bs{\alpha}}$ and induces asymmetry in 
$\ell^{\mathrm{D}}$. The 
parameter 
$\rho$ 
comes from the regular variation of $R$ at zero and infinity, respectively; this is reminiscent of 
the extremal attractors of 
elliptical distributions \cite{Opitz:2013}. The magnitude of $\rho$ has impact on the strength of 
dependence while its sign 
changes the overall shape of $\ell^{\mathrm{D}}$. Having $d+1$ parameters, the scaled extremal 
Dirichlet model may not be 
sufficiently rich to account for spatial dependence, unlike the H\"usler--Reiss or the extremal 
Student-$t$ models, which have 
one parameter for each pair of variables and are thus easily combined with distances. Also, it is 
less flexible than Dirichlet 
mixtures \cite{Boldi:2007}, which are however hard to estimate in high dimensions and require 
sophisticated machinery. 
To achieve greater flexibility, the scaled extremal Dirichlet model could perhaps 
be extended by working with more general scale mixtures, such as of the weighted Dirichlet 
distributions considered, e.g., in \cite{Hashorva:2008}.

Nonetheless, the scaled extremal Dirichlet model may naturally find applications whenever asymmetric 
extremal dependence is 
suspected; the latter may be caused, e.g., by causal relationships between the variables 
\cite{Genest/Neslehova:2013}.  
The 
stochastic structure of the scaled extremal Dirichlet model has several major advantages, that make 
the model easy to interpret, 
 estimate and simulate from. Its angular density has a simple form; in contrast to the asymmetric 
generalizations of the logistic 
and negative logistic models, this model does not place any mass on the vertices and 
lower-dimensional facets of the unit 
simplex. Another plus is the tractability of the de Haan representation and of the extremal functions, both 
expressible in terms of independent 
scaled Gamma variables; this allows for feasible inference  and stochastic simulation. While the 
scaled extremal Dirichlet stable 
tail dependence function $\ell^{\mathrm{D}}$ does not have a closed form in general, closed-form 
algebraic expressions exist  
when 
$\bs{\alpha}$ is integer-valued and in the bivariate case. Model selection for well-known families 
of extreme-value 
distributions can be performed through likelihood ratio tests. Another potentially useful feature is 
that $\rho \in (-\infty, 
\infty)$ can be allowed, with the convention that all variables whose indices $i$ are such that 
$-\rho  \le -\alpha_i$ are 
independent.

\subsection*{Acknowledgment}
Funding in partial support of this work was provided by the Natural Sciences and Engineering 
Research Council 
(RGPIN-2015-06801, CGSD3-459751-2014), 
the Canadian Statistical Sciences Institute, and the Fonds de recherche du Qu\'ebec -- Nature et 
technologies 
(2015--PR--183236).  We thank the acting Editor-in-Chief, Richard Lockhart, the Associate Editor and 
two anonymous 
referees for their valuable suggestions.

\setlength{\bibsep}{0.0pt}

\clearpage 
\begin{appendices}
\section{Proofs from Section \ref{sec:2}}\label{app:A}
{\noindent \it Proof of \Cref{prop:2}.}
To prove parts (a) and (b), recall that $1/X_i$ is distributed as $1/(RD_i)$, where $D_i \sim 
\mathsf{Beta}(\alpha_i, \bar 
\alpha   - \alpha_i)$ is independent of $R$. Furthermore, it is easy to show that $1/D_i \in 
\mathcal{M}(\Phi_{\alpha_i})$, 
which 
implies   that $\e(1/D_i^\beta) < \infty$ for any $\beta < \alpha_i$. The extremal behavior of 
$1/X_i$ will thus be determined 
by the   extremal behavior of either $1/R$ or $1/D_i$, depending on which one has a heavier tail. 
Indeed, 
Breiman's   Lemma \cite{Breiman:1965} 
implies that $1/X_i \in \mathcal{M}(\Phi_\rho)$ if $1/R \in 
\mathcal{M}(\Phi_\rho)$   for some $\rho < \alpha_i$ and that $1/X_i \in 
\mathcal{M}(\Phi_{\alpha_i})$ if $\e(1/R^{\alpha_i + 
\varepsilon}) < \infty$ for   some $\varepsilon > 0$. Finally, the fact that $1/X_i \in 
\mathcal{M}(\Phi_{\alpha_i})$ when $1/R 
\in   \mathcal{M}(\Phi_{\alpha_i})$ follows directly from the Corollary to Theorem~3 in 
\cite{Embrechts/Goldie:1980}. \qed

\medskip

The following lemma is a side result of \Cref{prop:2}, which is needed in the subsequent proofs.

\begin{lemma}\label{lem:svf} Suppose that $\bs{X} = R\bs{D}_\alpha$. If $1/R \in 
\mathcal{M}(\Phi_{\alpha_i})$ for some $i \in 
\{1, \dots, d\}$, then
\begin{align*}
\lim_{x \to \infty} \frac{\Pr\left({1}/{R} > x\right)} {\Pr\left({1}/{X_i} > x\right)}= 0.
\end{align*}
\end{lemma}

\medskip
{ \noindent \it Proof of \Cref{lem:svf}.}

Because $1/R \in \mathcal{M}(\Phi_{\alpha_i})$, $\Pr(1/R > x)$ is regularly varying with index 
$-\alpha_i$. In particular, for 
any 
$b \in (0,1)$, $\Pr(1/R>xb) / \Pr(1/R>x) \to b^{-\alpha_i}$ as $x \to \infty$. 
An application of Fatou's lemma thus gives
\begin{align*}
 \liminf_{x \to \infty} \frac{\Pr(1/X_i>x)}{\Pr(1/R>x)} & = 
 \liminf_{x \to \infty}  \int_0^1 \frac{\Pr(1/R>xb)}{\Pr(1/R>x)} f_{D_i}(b) \d b 
  \geq \int_0^1 
\frac{b^{-1} (1-b)^{\bar \alpha - \alpha_i - 1}}{\textrm{B}(\alpha_i, \bar \alpha - \alpha_i)}  \d 
b = \infty
\end{align*}
and hence the result. \qed

\section{Proofs from Section \ref{sec:3}}\label{app:B}
First recall the following property of the Dirichlet distribution, which is easily shown using the 
transformation 
formula for Lebesgue densities. 
\begin{lemma}\label{lem:B1} Let $ \bs{D}_\alpha$ be a Dirichlet random vector with parameters 
$\bs{\alpha}$. Then for any $2 
\le 
k \le d$ and any collection of distinct indices $1 \le i_1 < \dots < i_k \le d$, 
\begin{align*}
(D_{i_1}, \dots, D_{i_k}) \eqdis B_{i_1, \dots, i_k} \times \bs{D}_{(\alpha_{i_1}, \dots, 
\alpha_{i_k})}, 
\end{align*}
where $B_{i_1, \dots, i_k} \sim \mathsf{Beta}(\alpha_{i_1} + \cdots + \alpha_{i_k}, \bar \alpha - 
(\alpha_{i_1} + \cdots + 
\alpha_{i_k}))$ is independent of 
the $k$-variate Dirichlet vector $\bs{D}_{(\alpha_{i_1}, \dots, \alpha_{i_k})}$ with parameters 
$(\alpha_{i_1}, \dots, 
\alpha_{i_k})$. 
\end{lemma}

{ \noindent \it Proof of \Cref{thm:1}.} In order to prove part (a), recall that $\| \bs{X} \| 
\eqdis R$ is independent of $\bs{X} / 
\|\bs{X}\| \eqdis \bs{D}_{\bs{\alpha}}$.  Because $R \in \mathcal{M}(\Phi_\rho)$, there exists a 
sequence $(b_n)$ of constants 
in $(0, \infty)$ such that, for any Borel set $B \subseteq \Si_d$ and any $r > 0$, 
\begin{align*}
\lim_{n \to \infty} n\Pr\left(\|\bs{X}\| > b_nr, \frac{\bs{X}}{\|\bs{X}\|} \in B\right) = \lim_{n 
\to \infty} n 
\Pr(R>b_nr) \Pr(\bs{D}_{\bs{\alpha}} \in B)= r^{-\rho}\Pr(\bs{D}_{\bs{\alpha}} \in B).
\end{align*}
By Corollary 5.18 in \cite{Resnick:1987}, $\bs{X}\in \mathcal{M}(H_0)$ where for all $\bs{x} \in 
\mathbb{R}_+^d$, 
\begin{align*}
H_0(\bs{x})=\exp \left[- \e\left\{ \max\left( \frac{D_1^\rho}{x_1^\rho}, \dots, 
\frac{D_d^\rho}{x_d^\rho}\right) \right\} 
\right].
\end{align*}
Let $\Be(\cdot, \cdot)$ denote the Beta function. The univariate margins of $H_0$ are given, for 
all $i=1, \dots, 
d$ and $x > 0$, by
\begin{align*}
  F_{0i}(x) &=  \exp\left\{-x^{-\rho} \mathrm{E}(D_i^\rho)\right\} =\exp \left\{ -x^{-\rho}  
\frac{\mathrm{B}(\rho+\alpha_i, 
\bar 
\alpha-\alpha_i)}{\mathrm{B}(\alpha_i, 
\bar \alpha-\alpha_i)} \right\} = \exp \left\{ -x^{-\rho}  \frac{\Gamma(\alpha_i + \rho) 
\Gamma(\bar \alpha)}{\Gamma(\bar \alpha 
+ 
\rho) \Gamma(\alpha_i)} \right\}
\end{align*}
for $i\in \{1, \ldots, d\}$. The copula of $H_0$ then satisfies, for all $\bs{u} \in [0, 1]^d$, 
\begin{align*}
C_0(\bs{u}) = H_0\{ F_{01}^{-1}(u_1), \dots, F_{0d}^{-1}(u_d)\} = \exp\left(- \frac{\Gamma(\bar 
\alpha + \rho)}{\Gamma(\bar 
\alpha)} \e\left[ \max_{1 \le i \le d}\left\{ \frac{(-\log u_i) \Gamma(\alpha_i) 
D_i^\rho}{\Gamma(\alpha_i + \rho)}\right\} 
\right]\right).
\end{align*}
By \Cref{eq:evc}, the stable tail dependence function of $C_0$ thus indeed equals, for all $\bs{x} 
\in \mathbb{R}_+^d$, 
\begin{align*}
\ell(\bs{x}) = \frac{\Gamma(\bar \alpha + \rho)}{\Gamma(\bar \alpha)} \e\left[ \max\left\{ 
\frac{x_1 \Gamma(\alpha_1) 
D_1^\rho}{\Gamma(\alpha_1 + \rho)}, \dots, \frac{x_d \Gamma(\alpha_d) D_d^\rho}{\Gamma(\alpha_d + 
\rho)}\right\} \right].
\end{align*}
\smallskip
The part (b) follows directly from Proposition~2.2 in \cite{Hashorva:2015} upon setting $p=1$ and 
taking, for 
$i=1,\dots, d$ and $j=1,\dots, d$, $\lambda_{ij}=1$ whenever $i=j$ and $\lambda_{ij}=0$ otherwise.

\smallskip
To prove part (c), recall first that from 
\Cref{prop:1}, for $i=1, \dots, d$, $X_i \in \mathcal{M}(\Psi_{\rho + \bar \alpha -\alpha_i})$ and 
hence there 
exist sequences 
$(a_{ni})\in (0, \infty)$, 
$(b_{ni}) \in \mathbb{R}$, such that for all $x \in \mathbb{R}$, 
\begin{align*}
\lim_{n\to \infty} n \Pr(X_i > a_{ni} x + b_{ni}) = -\log \{\Psi_{\rho + \bar \alpha 
-\alpha_i}(x)\}.                             
\end{align*}
Next, observe that as in the proof of Proposition 5.27 in \cite{Resnick:1987}, $\bs{X} \in 
\mathcal{M}(H_0)$ follows if for all 
$1 \le i < j \le d$ and $x_k$ such that $\Psi_{\rho + \bar \alpha -\alpha_k} (x_k) >0$ for $k=i, 
j$, 
\begin{equation}\label{eq:BR1}
\lim_{n\to \infty} n \Pr(X_i > a_{ni} x_i + b_{ni}, X_j > a_{nj} x_j + b_{nj}) =0.
\tag{B.1}
\end{equation}
To prove that \eqref{eq:BR1} indeed holds, it suffices to assume that $d=2$. This is because for 
arbitrary indices $1 \le i < j 
\le d$, \Cref{lem:B1} implies that $(X_i, X_j) \eqdis R^*(B, 1-B)$, where $B \sim 
\mathrm{Beta}(\alpha_i, \alpha_j)$, $R^* 
\eqdis 
R Y$ is independent of $B$ and $Y \sim \mathrm{Beta}(\alpha_i + \alpha_j, \bar \alpha - \alpha_i 
-\alpha_j)$ is independent of 
$B$ and $R$. Because $\Pr(R^* \le 0)=0$, Theorem~4.5 in \cite{Hashorva/Pakes:2010} implies that 
$R^* \in \mathcal{M}(\Psi_{\rho + 
\bar \alpha - \alpha_i - \alpha_j})$ when $R \in \mathcal{M}(\Psi_{\rho})$ for some $\rho >0$. 
Thus suppose that $d=2$ and write $(D_1, D_2)\equiv (B, 1-B)$, where $B \sim 
\mathsf{Beta}(\alpha_1, \alpha_2)$. Fix 
arbitrary 
$x_1, x_2 \in \mathbb{R}$ are such that $\Psi_{\rho + \bar \alpha -\alpha_i} (x_i) >0$ for $i=1, 
2$. Then because for any $a, c 
>0$ and $b \in (0, 1)$, 
$\max\{a/b, c/(1-b) \} \ge a+c$, one has
\begin{align*}
0 \le \Pr(X_1 > a_{n1} x_1 + b_{n1}, X_2  > a_{n2} x_2 + b_{n2}) &= \Pr\left\{ R > \max 
\left(\frac{a_{n1} x_1  +b_{n1}}{B}, 
\frac{a_{n2} x_2  +b_{n2}}{1-B}\right)\right\} \\
&\le \Pr(R > a_{n1} x_1 + b_{n1} + a_{n2} x_2 + b_{n2}).
\end{align*}
In order to prove \Cref{eq:BR1}, it thus suffices to show that 
\begin{equation}\label{eq:BR3}
\lim_{n\to\infty}n\Pr(R > a_{n1} x_1 + b_{n1} + a_{n2} x_2 + b_{n2}) = 0.
\tag{B.2}
\end{equation} 
This however follows immediately from the fact that if $R \in \mathcal{M}(\Psi_\rho)$ for some 
$\rho > 0$, the upper end-point 
$r$ 
of $R$, viz. $r= \sup\{x : \Pr(R \le x) < 1\}$, is finite. Because for $i=1, 2, $ $r$ is also the 
upper 
endpoint of $X_i$, $a_{ni} x_i + b_{ni} \to r$ as $n\to \infty$. This means that there exists $n_0 
\in \mathbb{N}$ so that 
for all $n \ge n_0$, $a_{n1} x_1 + b_{n1} + a_{n2} x_2 + b_{n2} > r$ and $\Pr(R > a_{n1} x_1 + 
b_{n1} + a_{n2} x_2 + b_{n2})=0$. 
This proves 
\Cref{eq:BR1} and hence also Theorem \ref{thm:1} (c). Note that alternatively, part (c) could be  
proved using Theorem~2.1 
in \cite{Hashorva:2015} similarly to the proof of Proposition~2.2 therein.
\qed

\bigskip
The proof of Theorem \ref{thm:2} requires the following  technical lemma. 
\begin{lemma}\label{lem:B2} Suppose that $\bs{D}_{\bs{\alpha}} =(D_{1}, \dots, D_{d})$ is a 
Dirichlet random vector with 
parameters $\bs{\alpha}$. Further let  $R$ be a positive random variable independent of 
$\bs{D}_{\bs{\alpha}}$ such that $\Pr(R 
\le 0) =0$, and let $\bs{X} = R \bs{D}_\alpha$. Then for any $1 \le i < j \le d$ and any $x_i, x_j 
\in (0, \infty)$, 
\begin{align*}
\lim_{n\to \infty} n \Pr\left( \frac{1}{X_i} > a_{ni} x_i, \frac{1}{X_j} > a_{nj} x_j\right) =0
\end{align*}
if either:
\begin{enumerate}
\item[\textrm{(i)}] $1/R \in \mathcal{M}(\Phi_\rho)$ with $\rho \in [\alpha_i \wedge \alpha_j, 
\alpha_1\vee\alpha_2]$, and for 
$k=i, j$, $(a_{nk})$ is a sequence of positive constants such that $n \Pr(1/X_k > a_{nk} x_k) \to 
x_k^{-(\alpha_k \wedge \rho)}$ 
as $n\to \infty$;
\item[\textrm{(ii)}] $\textrm{E}(1/R^\beta) < \infty$ for some $\beta > \alpha_i \vee\alpha_j$  and 
for $k=i, j$, $(a_{nk})$ is a 
 sequence of positive constants such that $n \Pr(1/X_k > a_{nk} x_k)  \to x_k^{-\alpha_k}$ as $n\to 
\infty$.
\end{enumerate}
\end{lemma}
{ \noindent \it Proof of Lemma \ref{lem:B2}.}  Observe first that when $d >2$, Lemma~\ref{lem:B1} 
implies that $(X_{i}, X_{j})\eqdis 
R^*(B, 1-B)$, where $R^* \indep B$, $B \sim \mathsf{Beta}(\alpha_i, \alpha_j)$ and $R^* = RY$, with 
$Y \indep R$ and $Y \sim 
\mathsf{Beta}(\alpha_i + \alpha_j, \bar \alpha - \alpha_i- \alpha_j)$. Now note that $1/Y \in 
\mathcal{M}(\Phi_{\alpha_i + 
\alpha_j})$. Thus if $1/R \in \mathcal{M}(\Phi_\rho)$ for $\rho \in [\alpha_i \wedge \alpha_j, 
\alpha_1\vee\alpha_2]$, $\rho < 
\alpha_i + \alpha_j$ and Breiman's Lemma implies that $1/R^* \in \mathcal{M}(\Phi_\rho)$. 
Further, if $\textrm{E}(1/R^\beta) < \infty$ 
for some $\beta \in (\alpha_i \vee\alpha_j, \alpha_i + \alpha_j)$, $\textrm{E}\{1/(R^*)^\beta \} < 
\infty$ given that 
$\textrm{E}(1/Y^\beta)
<\infty$. We can thus assume without loss of generality that $d=2$ and $\alpha_1 \le \alpha_2$; we 
shall also write 
$(D_1, D_2)\equiv (B, 1-B)$, where $B \sim \mathsf{Beta}(\alpha_1, \alpha_2)$.

\medskip
To prove part (i), note first that the existence of the sequences $(a_{nk})$, $k=i, j$, follows 
from \Cref{prop:2}, by which 
$1/X_k \in \mathcal{M}(\Phi_{\rho \wedge \alpha_k})$ for $k=i, j$, and the Poisson approximation 
\citep[Proposition 
3.1.1]{Embrechts/Kluppelberg/Mikosch:1997}. Next, observe that for any constants $a, c > 0$ and any 
$b \in (0, 1)$, 
\begin{equation}\label{eq:B1}
\frac{ac}{a+c} \le ab \vee c(1-b) < a \vee c.
\tag{B.3}
\end{equation}
Indeed, when $b < c/(a+c)$, $ab \vee c(1-b) = c(1-b)$ and $c(1-b) \in (ac/(a+c), c)$, while when $b 
\ge c/(a+c)$, $ab \vee 
c(1-b) 
= ab$ and $ab \in [ac/(a+c), a)$. To show the claim in part (i), distinguish the cases below:

\smallskip
{\it Case I.} $\alpha_1 = \alpha_2$. Here, $\rho = \alpha_1 =\alpha_2$ and $X_1 \eqdis X_2$, so 
that $a_{n1}/a_{n2} \to 1$ by 
the 
Convergence to Types Theorem \cite{Resnick:1987}. By \Cref{eq:B1}, 
\begin{equation}\label{eqBL1}
0 \le n \Pr\left( \frac{1}{X_1} > a_{n1} x_1, \frac{1}{X_2} > a_{n2} x_2\right) = n \Pr\left[ 
\frac{1}{R} > \max\{a_{n1} x_1 B, 
a_{n2} x_2(1-B)\}\right]  \le n \Pr\left( \frac{1}{R} >  \frac{a_{n1} a_{n2} x_1 x_2}{a_{n1} x_1 + 
a_{n2} x_2} \right).
\tag{B.4}
\end{equation}  
Because $(x_1 x_2)/\{(a_{n1}/a_{n2} ) x_1 + x_2\} \to (x_1x_2)/(x_1+x_2)$ as $n\to \infty$, 
\begin{align*}
\lim_{n\to \infty} n \Pr\left( \frac{1}{X_1} >  \frac{a_{n1} a_{n2} x_1 x_2}{a_{n1} x_1 + a_{n2} 
x_2} \right) = \left(\frac{x_1 
x_2}{x_1 + x_2}\right)^{- \alpha_1}.
\end{align*}
Furthermore, by \Cref{lem:svf}, given that $ (a_{n1} a_{n2} x_1 x_2)/(a_{n1} x_1 + a_{n2} x_2) \to 
\infty$ as $n\to \infty$, 
\begin{align*}
\lim_{n\to\infty}\dfrac{\Pr\left( {1}/{R} >  \frac{a_{n1} a_{n2} x_1 x_2}{a_{n1} x_1 + a_{n2} x_2} 
\right)}{\Pr\left( 
{1}/{X_1} >  \frac{a_{n1} a_{n2} x_1 x_2}{a_{n1} x_1 + a_{n2} x_2} \right)} =0, 
\end{align*}
so that the right-hand side in \Cref{eqBL1} tends to $0$ as $n\to \infty$, and this implies the 
claim.

\smallskip
{\it Case II.} $\alpha_1 < \alpha_2$ and $\rho = \alpha_2$. Then for $i=1, 2$, there exists a 
slowly varying function $L_i$ such 
that $a_{ni} = n^{1/\alpha_i} L_i(n)$. Hence $a_{n2}/a_{n1} \to 0$ and  $(x_1 x_2)/\{ x_1 + x_2 
(a_{n2}/a_{n1} )\} \to x_2$ as 
$n\to \infty$. Consequently, 
\begin{align*}
\lim_{n\to \infty} n \Pr\left( \frac{1}{X_2} >  \frac{a_{n1} a_{n2} x_1 x_2}{a_{n1} x_1 + a_{n2} 
x_2} \right) = x_2^{- 
\alpha_2}.
\end{align*}
Moreover, by \Cref{lem:svf}, given that $ (a_{n1} a_{n2} x_1 x_2)/(a_{n1} x_1 + a_{n2} x_2) \to 
\infty$ as $n\to \infty$, 
\begin{align*}
\lim_{n\to\infty}\dfrac{\Pr\left( 1/R >  \frac{a_{n1} a_{n2} x_1 x_2}{a_{n1} x_1 + a_{n2} x_2} 
\right)}{\Pr\left( 
1/{X_2} >  \frac{a_{n1} a_{n2} x_1 x_2}{a_{n1} x_1 + a_{n2} x_2} \right)} =0, 
\end{align*}
so that again the right-hand side in \Cref{eqBL1} tends to $0$ as $n\to \infty$.

\smallskip
{\it Case III.} $\alpha_1 < \alpha_2$ and $\rho \in [\alpha_1, \alpha_2)$. In this case, $1/X_1 \in 
\mathcal{M}(\Phi_{\alpha_1})$ 
and $1/X_2 \in \mathcal{M}(\Phi_\rho)$. Therefore, either directly when $\rho > \alpha_1$ or by 
\Cref{lem:svf}, one can easily 
deduce that
\begin{align}\label{eq:BL4}
\lim_{x\to \infty} \frac{\Pr(1/R > x)}{\Pr(1/X_1 > x)} = 0.
\tag{B.5}
\end{align}
At the same time, Breiman's Lemma \cite{Breiman:1965} implies that 
\begin{align}\label{eq:BL5}
\lim_{x\to \infty} \frac{\Pr(1/R > x)}{\Pr(1/X_2 > x)} = \frac{1}{\mathrm{E}\{1/(1-B)^\rho\}} 
=\frac{\mathrm{B}(\alpha_1, \alpha_2)}{\mathrm{B}(\alpha_1, \alpha_2-\rho)}.
\tag{B.6}
\end{align}
Hence, for any $b \in (0, 1)$, the limit of $n \Pr\left\{1/R > a_{n2} x_2 (1-b)\right\}$ as $n\to 
\infty$ equals
\begin{align*}
 \lim_{n\to \infty}  n \Pr\left\{1/X_2 > a_{n2} x_2 
(1-b)\right\} 
\frac{\Pr\left\{1/R > a_{n2} x_2 (1-b)\right\}}{\Pr\left\{1/X_2 > a_{n2} x_2 (1-b)\right\}} 
=\{x_2(1-b)\}^{-\rho} 
\frac{\mathrm{B}(\alpha_1, \alpha_2)}{\mathrm{B}(\alpha_1, \alpha_2-\rho)}
\end{align*}
so that
\begin{multline}\label{eq:BL2}
\lim_{n\to \infty} \int_0^1  n \Pr\{1/R > a_{n2} x_2 (1-b)\} \frac{b^{\alpha_1 
-1}(1-b)^{\alpha_2-1}}{\textrm{
B}(\alpha_1, \alpha_2)} 
\d b = n \Pr(1/X_2 > a_{n2} x_2) = x_2^{-\rho}\\ = x_2^{-\rho} \int_0^1  \frac{b^{\alpha_1 
-1}(1-b)^{\alpha_2-\rho-1}}{\textrm{
B}(\alpha_1, \alpha_2-\rho)} \d b = \int_0^1 \lim_{n\to \infty} n \Pr\{1/R > a_{n2} x_2 (1-b)\} 
\frac{b^{\alpha_1 
-1}(1-b)^{\alpha_2-1}}{\textrm{B}(\alpha_1, \alpha_2)} \d b.
\tag{B.7}
\end{multline}
Given that for any $b \in (0, 1)$, $\Pr\{1/R > a_{n1} x_1 b, 1/R > a_{n2}x_2 (1-b)\} \le \Pr\{1/R > 
a_{n2}x_2 (1-b)\}$, 
\begin{align*}
\lefteqn{\int_0^1 \liminf_{n\to \infty} \bigl( n \bigl[ \Pr\{1/R > a_{n2} x_2(1-b)\} - \Pr\{1/R > 
a_{n1}x_1 b, 1/R > a_{n2}x_2 
(1-b)\}\bigr]\bigr)  \frac{b^{\alpha_1 -1}(1-b)^{\alpha_2-1}}{\textrm{B}(\alpha_1, \alpha_2)}\d b} 
& \\
&\le \liminf_{n\to \infty} \int_0^1  n \bigl[ \Pr\{1/R > a_{n2} x_2(1-b)\} - \Pr\{1/R > a_{n1} 
x_1b, 1/R > a_{n2} 
x_2(1-b)\}\bigr] \frac{b^{\alpha_1 -1}(1-b)^{\alpha_2-1}}{\textrm{B}(\alpha_1, \alpha_2)} \d b
\end{align*}
by Fatou's Lemma. Because of \Cref{eq:BL2}, this inequality simplifies to
\begin{align*}
\lefteqn{x_2^{-\rho} - \int_0^1 \limsup_{n\to \infty}\bigl[n \Pr\{1/R > a_{n1}x_1 b, 1/R > a_{n2} 
x_2(1-b)\} \bigr] 
\frac{b^{\alpha_1 -1}(1-b)^{\alpha_2-1}}{\textrm{B}(\alpha_1, \alpha_2)} \d b} & \\
&\le x_2^{-\rho} - \limsup_{n\to \infty} \int_0^1 n \Pr
\{1/R > a_{n1} x_1b, 1/R > a_{n2}x_2 (1-b)\} \frac{b^{\alpha_1 
-1}(1-b)^{\alpha_2-1}}{\textrm{B}(\alpha_1, \alpha_2)} \d b
\end{align*}
and hence 
\begin{multline*}
0 \le \limsup_{n\to \infty}\bigl\{n \Pr(1/X_1 > a_{n1} x_1, 1/X_2 > a_{n2}x_2)\bigr\} \\\le 
\int_0^1 \limsup_{n\to 
\infty}\bigl[n \Pr\{1/R > a_{n1}x_1 b, 1/R > a_{n2}x_2 (1-b)\} \bigr]\frac{b^{\alpha_1 
-1}(1-b)^{\alpha_2-1}}{\textrm{
B}(\alpha_1, \alpha_2)} \d b.
\end{multline*}
To show the desired claim, it thus suffices to show that for arbitrary $b \in (0, 1)$, 
\begin{align}\label{eq:BL3}
\lim_{n \to \infty} n \Pr\{1/R > a_{n1} x_1b, 1/R > a_{n2} x_2(1-b)\}  = 0.
\tag{B.8}
\end{align}
To this end, fix $b \in (0, 1)$ and observe that $a_{n1}/a_{n2} \to \infty$. Indeed, if $\rho > 
\alpha_1$, this follows 
directly 
from the fact that $a_{n1} = n^{1/\alpha_1} L_1(n)$ and $a_{n2} = n^{1/\rho} L_2(n)$ for some 
slowly varying functions $L_1, 
L_2$. When $\rho = \alpha_1$, suppose that $\liminf_{n\to \infty} a_{n1} /a_{n2}$ were finite. Then 
there exists a subsequence 
$a_{n_k 1}/a_{n_k 2}$ such that $a_{n_k 1}/a_{n_k 2} \to a$ as $k\to \infty$ for some $a \in [0, 
\infty)$. Hence, for a fixed 
$\varepsilon > 0$ and all $k \ge k_0$, $a_{n_k 1}/a_{n_k 2} \le a+\varepsilon$. Using the latter 
observation and \Cref{eq:BL5}, 
\begin{multline*}
\lim_{k \to \infty} n_k \Pr(1/R > a_{n_k 1}) \ge  \lim_{k \to \infty} n_k \Pr\{1/R > a_{n_k 
2}(a+\varepsilon)\} \\  = \lim_{k 
\to 
\infty} n_k \Pr\{1/X_2 > a_{n_k 2}(a+\varepsilon)\} \frac{\Pr\{1/R > a_{n_k 
2}(a+\varepsilon)\}}{\Pr\{1/X_2 > a_{n_k 
2}(a+\varepsilon)\}} = (a+\varepsilon)^{-\rho} \frac{\mathrm{B}(\alpha_1, 
\alpha_2)}{\mathrm{B}(\alpha_1, \alpha_2-\rho)} > 0.
\end{multline*}
At the same time, by \Cref{eq:BL4}, 
\begin{align*}
\lim_{k \to \infty} n_k \Pr(1/R > a_{n_k 1}) = \lim_{k \to \infty} n_k \Pr(1/X_1 > a_{n_k 1}) 
\frac{\Pr(1/R > a_{n_k 
1})}{\Pr(1/X_1 > a_{n_k 1})} =0
\end{align*}
and hence a contradiction. Therefore, $\liminf_{n\to \infty} a_{n1} /a_{n2} =\infty$ and hence 
$a_{n1} /a_{n2} \to 
\infty$ 
as $n \to \infty$.
Because $a_{n1} b > a_{n2} (1-b)$ if and only if $b > a_{n2}/(a_{n1} + a_{n2})$ and $a_{n2}/(a_{n1} 
+ a_{n2}) \to 0$ as $n \to 
\infty$, there exists $n_0$ such that for all $n \ge n_0$, 
\begin{align*}
n \Pr\{1/R > a_{n1} x_1b, 1/R > a_{n2} x_2(1-b)\} = n \Pr(1/R > a_{n1} x_1 b)  = n \Pr(1/X_1 > 
a_{n1} x_1 b) \frac{\Pr(1/R > 
a_{n1} x_1 b)}{\Pr(1/X_1 > a_{n1} x_1 b)}.
\end{align*}
The last expression tends to $0$ as $n \to \infty$ by \Cref{eq:BL4} and hence \Cref{eq:BL3} indeed 
holds.

\medskip
To prove part (ii), first recall that by \Cref{prop:2}~(b), $1/X_i \in 
\mathcal{M}(\Phi_{\alpha_i})$, $i=1, 2$, and hence the scaling sequences $(a_{n1})$ and $(a_{n2})$ indeed exist. Recall that for 
$i=1, 2$, $a_{ni} = 
n^{1/\alpha_i} L_i(n)$ for some slowly 
varying function $L_i$.
As in the proof of part (i), $n\Pr(1/X_1 > a_{n1} x_1, 1/X_2 > a_{n2} x_2)$ can be bounded above by 
the right-hand side in 
\Cref{eqBL1}. Markov's inequality further implies that for $\beta \in (\alpha_2, \alpha_1 + 
\alpha_2)$ such that 
$\textrm{E}(1/R^\beta)< \infty$, 
\begin{align*}
n \Pr\left( \frac{1}{R} > \frac{a_{n1} a_{n2} x_1 x_2}{a_{n1} x_1 + a_{n2} x_2} \right) & \le n \E{ 
1/R^{\beta}} \frac{(a_{n1} 
x_1 + a_{n2} x_2)^{\beta}}{(a_{n1} a_{n2} x_1 x_2)^{\beta}} =  
\frac{\E{1/R^\beta}}{(x_1x_2)^{\beta}} 
\left\{ \frac{x_1}{n^{1/\alpha_2 -1/\beta}L_2(n)} + \frac{x_2}{n^{1/\alpha_1 - 1/\beta} 
L_1(n)}\right\}^{\beta}
\end{align*}
The right-most expression tends to $0$ as $n\to \infty$ because for any $i=1, 2$ and $\rho > 0$, 
$n^\rho L_i(n) \to \infty$.
\qed

\bigskip
\noindent
{ \noindent \it Proof of Theorem \ref{thm:2}.}  
First note that a positive random vector $\bs{Y}$ is in the maximum domain of attraction of a 
multivariate extreme-value 
distribution $H_0$ with Fr\'echet margins if and only if there exist sequences of positive 
constants $(a_{ni}) \in (0, \infty)$, 
$i=1, \dots, d$, 
so that, for all $\bs{y} \in \mathbb{R}_+^d$, 
\begin{multline}\label{eq:B2}
\lim_{n\to \infty} n \{ 1- \Pr(Y_1 \le a_{n1} y_1, \dots, Y_d \le a_{nd} y_d)\} = \\
\lim_{n\to \infty} n \left\{ \sum_{k=1}^d \sum_{1 \le i_1 < \dots <  i_k\le d} 
(-1)^{k+1}\Pr(Y_{i_1} > a_{ni_1} y_{i_1}, \dots, 
Y_{i_k} > a_{n i_k} y_{i_k})\right\} = - \log H_0(\bs{y}).
\tag{B.9}
\end{multline}
This multivariate version of the Poisson approximation holds by the same argument as in the 
univariate case \citep[Proposition 
3.1.1]{Embrechts/Kluppelberg/Mikosch:1997}.

To prove part (a), suppose that $1/R \in \mathcal{M}(\Phi_\rho)$ for some $\rho \in (0, \alpha_M]$. 
By \Cref{prop:2}, one then 
has that for any $i\in \mathbb{I}_1$, $1/(RD_i) \in \mathcal{M}(\Phi_{\alpha_i})$. For any $i \in 
\mathbb{I}_1$, let $(a_{ni})$ 
be 
a sequence of positive constants such that, for all $x > 0$, $n \Pr\{1/(RD_i) > a_{ni} x\} \to 
x^{-\alpha_i}$ as $n \to \infty$; 
such a sequence exists by the univariate Poisson approximation \citep[Proposition 
3.1.1]{Embrechts/Kluppelberg/Mikosch:1997}. 
The 
same result also guarantees the existence of a sequence $(a_n)$ of positive constants such that, 
for all $x > 0$, $n \Pr(1/R > 
a_n 
x) \to x^{-\rho}$ as $n \to \infty$. Now set, for any $i \in \mathbb{I}_2$, 
\begin{align}
b_i = \E{D_{i}^{-\rho}} = \frac{\Gamma(\alpha_i - \rho) \Gamma(\bar \alpha - \alpha_i)}{\Gamma(\bar 
\alpha - \rho)} \times 
\frac{\Gamma(\bar \alpha)}{\Gamma(\alpha_i) \Gamma(\bar \alpha - \alpha_i)} = \frac{\Gamma(\bar 
\alpha) / \Gamma(\bar \alpha - 
\rho)}{\Gamma(\alpha_i) / \Gamma(\alpha_i - \rho)}, 
\label{eq:tB11}\tag{B.10}
\end{align}
and define, for any $i \in \mathbb{I}_2$ and $n \in \mathbb{N}$, $a_{ni} =  b_i^{1/\rho} a_n$. As 
detailed in the proof of 
\Cref{prop:2}~(a), Breiman's  Lemma then implies that, for all $i\in \mathbb{I}_2$ and $x > 0$, 
\begin{align*}
\lim_{n\to \infty} n \Pr\left\{ \frac{1}{RD_{i}} > a_{ni} x\right\}= \lim_{n\to \infty} n  
\Pr\left\{ \frac{1}{R} > a_n 
(b_i^{1/\rho}x)\right\} \frac{\Pr\left\{ \frac{1}{RD_{i}} > a_n (b_i^{1/\rho}x)\right\}}{ 
\Pr\left\{ \frac{1}{R} > a_n 
(b_i^{1/\rho}x)\right\}} =x^{-\rho} b_i^{-1} b_i = x^{-\rho}, 
\end{align*}
given that for all $i\in \mathbb{I}_2$, $D_{i} \sim \mathsf{Beta}(\alpha_i, \bar \alpha - 
\alpha_i)$. 

Next, fix an arbitrary $\bs{x} \in (0, \infty)^d$, $k \in \{ 2, \dots, d\}$ and indices $1 \le i_1 
< \dots < i_k \le d$. To 
calculate the limit of $n \Pr ( 1/(R D_{i_1}) > a_{ni_1} x_{i_1}, \dots, 1/(R D_{i_k}) > a_{ni_k} 
x_{i_k})$, two cases must 
be distinguished: 

\smallskip

{\it Case I.} $\{i_1, \dots, i_k\} \cap \mathbb{I}_1 \neq \emptyset$. In this case, suppose, 
without loss of generality, that 
$i_1 
\in \mathbb{I}_1$. Then
\begin{align*}
0\le n \Pr \left( \frac{1}{R D_{{i_1}}} > a_{ni_1} x_{i_1}, \dots, \frac{1}{R D_{{i_k}}} > a_{ni_k} 
x_{i_k}\right) \le n \Pr 
\left( \frac{1}{R D_{{i_1}}} > a_{ni_1} x_{i_1}, \frac{1}{R D_{{i_2}}} > a_{ni_2} x_{i_2}\right). 
\end{align*}
Now either $i_2 \in \mathbb{I}_1$, in which case $\rho \ge \alpha_{i_1} \vee \alpha_{i_2}$, or $i_2 
\in \mathbb{I}_2$, so that 
$\alpha_{i_1} \le \rho < \alpha_{i_2}$. Either way, Lemma \ref{lem:B2} implies that 
\begin{align*}
\lim_{n\to \infty} n \Pr \left( \frac{1}{R D_{{i_1}}} > a_{ni_1} x_{i_1}, \frac{1}{R D_{{i_2}}} > 
a_{ni_2} x_{i_2}\right) =0
\end{align*}
and consequently $n \Pr \{ 1/(R D_{i_1}) > a_{ni_1} x_{i_1}, \dots, 1/(R D_{i_k}) > a_{ni_k} 
x_{i_k}\} \to 0$ as $n\to \infty$.

{\it Case II.} $\{i_1, \dots, i_k\} \cap \mathbb{I}_1 = \emptyset$. In this case, let $Z_{i_1, 
\dots, i_k} = \max(x_{i_1} 
(b_{i_1})^{1/\rho} D_{{i_1}}, \dots, x_{i_k} (b_{i_k})^{1/\rho} D_{{i_k}})$ and observe that for 
any $\varepsilon > 0$ such 
that 
$\rho + \varepsilon < \min(\alpha_1, \dots, \alpha_d)$, 
\begin{align*}
\e\left( \frac{1}{Z_{i_1, \dots, i_k}^{\rho + \varepsilon}} \right) \le x_{i_1}^{-\rho - 
\varepsilon} b_{i_1}^{-(\rho + 
\varepsilon)/\rho} \e\left( \frac{1}{D_{{i_1}}^{\rho + \varepsilon}} \right) < \infty.
\end{align*}
Therefore, by Breiman's Lemma, 
\begin{multline*}
\lim_{n \to \infty} n \Pr \left( \frac{1}{R D_{{i_1}}} > a_{ni_1} x_{i_1}, \dots, \frac{1}{R 
D_{{i_k}}} > a_{ni_k} 
x_{i_k}\right)  = \lim_{n \to \infty} n \Pr \left( \frac{1}{R Z_{i_1, \dots, i_k} } > a_n \right) = 
\E{ Z_{i_1, \dots, 
i_k}^{-\rho}} 
\\
= {\mathrm E}\left[ \left\{  \max_{1 \le j \le k}\bigl(x_{i_j} b_{i_j}^{1/\rho} D_{{i_j}}\bigr) 
\right\}^{-\rho} \right] 
=  
{\mathrm E}\left[  \min_{1 \le j \le k}\left\{\frac{\left(x_{i_j} 
D_{{i_j}}\right)^{-\rho}}{b_{i_j}}\right\} \right].
\end{multline*}
Putting the above calculations together, one then has, for any $\bs{x} \in \mathbb{R}_+^d$, 
\begin{align*}
\lim_{n \to \infty} n\left\{ 1-  \Pr \left( \frac{1}{R D_{{1}}} \le a_{n1} x_{1}, \dots, \frac{1}{R 
D_{{d}}} \le a_{nd} 
x_{d}\right) \right\} =\sum_{i \in \mathbb{I}_1} x_{i}^{-\alpha_i} + \sum_{k=1}^{|\mathbb{I}_2|} 
\sum_{\substack{\{i_1, \dots, i_k\} \subseteq \mathbb{I}_2 \\ i_1 < \dots <  i_k}} (-1)^{k+1} 
{\mathrm E}\left[ 
\min_{1 \le j \le k}\left\{\frac{\left(x_{i_j}D_{{i_j}}\right)^{-\rho}}{b_{i_j}} \right\} \right]
\end{align*}
Furthermore, one can readily establish by induction that for any $\bs{t} \in \mathbb{R}^d$, 
\begin{align*}
\sum_{k=1}^{|\mathbb{I}_2|} \sum_{\substack{\{i_1, \dots, i_k\} \subseteq \mathbb{I}_2 \\ i_1 < 
\dots <  i_k}} (-1)^{k+1} 
\min(t_{i_1}, \dots, t_{i_k}) = \max_{i \in \mathbb{I}_2}(t_i).
\end{align*}
Hence, for any $\bs{x} \in \mathbb{R}_+^d$, 
\begin{align*}
\lim_{n \to \infty} n\left\{ 1-  \Pr \left( \frac{1}{R D_{{1}}} \le a_{n1} x_{1}, \dots, \frac{1}{R 
D_{{d}}} \le a_{nd} 
x_{d}\right) \right\}  =  \sum_{i \in \mathbb{I}_1} x_{i}^{-\alpha_i} + {\mathrm E} \left[\max_{i 
\in 
\mathbb{I}_2}\left\{\frac{\left(x_{i}D_i\right)^{-\rho}}{b_i} \right\}\right].
\end{align*}
By the multivariate Poisson approximation \eqref{eq:B2}, $1/\bs{X} \in \mathcal{M}(H_0)$, 
where for all $\bs{x} \in \mathbb{R}_+^d$, 
\begin{align*}
H_0(\bs{x}) = \exp \left( - \sum_{i \in \mathbb{I}_1} x_{i}^{-\alpha_i} - {\mathrm E} \left[\max_{i 
\in 
\mathbb{I}_2}\left(\frac{(x_{i}D_i)^{-\rho}}{b_i} \right\}\right] \right).
\end{align*}
The univariate margins of $H_0$ are given, for all $i\in \mathbb{I}_1$, by $F_{0i}(x) = 
x^{-\alpha_i}$ and for all 
$i \in \mathbb{I}_2$, 
$F_{0i}(x) = \exp(-x^{-\rho})$. By Sklar's Theorem, the unique copula of $H_0$ is given, for all 
$\bs{u} \in [0, 1]^d$, by 
\eqref{eq:evc}, where for all $\bs{x} \in \mathbb{R}_+^d$, 
\begin{align*}
\ell(\bs{x}) = \sum_{i\in \mathbb{I}_1} x_i + {\mathrm E} \left\{\max_{i \in \mathbb{I}_2}\left( 
\frac{x_iD_{i}^{-\rho}}{b_i} 
\right)\right\}.
\end{align*}
The first expression for $\ell$ follows immediately from \Cref{eq:tB11}. The second 
expression is readily verified using  \Cref{lem:B1}, given the fact that  if $B \sim 
\textrm{Beta}(\bar \alpha_2, \bar \alpha - 
\bar 
\alpha_2)$, $\textrm{E}(B^{-\rho}) = \Gamma(\bar \alpha_2 - \rho) \Gamma(\bar \alpha)/\Gamma(\bar 
\alpha - \rho) \Gamma(\bar 
\alpha_2)$.

\medskip
To prove part (b), recall that by Proposition \ref{prop:2}~(b), $1/X_i \in 
\mathcal{M}(\Phi_{\alpha_i})$. Hence, there exist 
sequences of positive constants $(a_{ni})$, $i=1, \dots, d$, such that for all $i=1, \dots, d$ and 
all $x > 0$, $n 
\Pr(1/(RD_i) > 
a_{ni} x) \to x^{-\alpha_i}$ as $n \to \infty$. By Lemma \ref{lem:B1}~(ii), it also follows that 
for  arbitrary $\bs{x} \in 
(0, \infty)^d$, $k \in \{ 2, \dots, d\}$ and indices $1 \le i_1 < \dots < i_k \le d$, 
\begin{align*}
0\le \lim_{n\to \infty} n \Pr \left( \frac{1}{R D_{{i_1}}} > a_{ni_1} x_{i_1}, \dots, \frac{1}{R 
D_{{i_k}}} > a_{ni_k} 
x_{i_k}\right)  \le \lim_{n\to \infty} n \Pr \left( \frac{1}{R D_{{i_1}}} > a_{ni_1} x_{i_1}, 
\frac{1}{R D_{{i_2}}} > 
a_{ni_2} 
x_{i_2}\right) = 0.
\end{align*}
Thus, by \Cref{eq:B2}, $1/\bs{X}$ is in the domain of attraction of the multivariate extreme-value 
distribution given, for all 
$\bs{x} \in \mathbb{R}_+^d$, by $H_0(\bs{x}) = \exp(-x_1^{-\alpha_1} - \dots - x_d^{-\alpha_d})$, 
as was to be showed.
\qed

\section{Proofs from Section \ref{sec:4}}\label{app:C}

{ \noindent \it Proof of \Cref{prop:3}.} In view of \Cref{cor:1} and \Cref{thm:2} in 
\cite{Larsson/Neslehova:2011}, it only remains to 
derive the explicit expression for $\ell^{\mathrm{nD}}$. Because $1-\psi(1/\cdot) \in 
\mathcal{R}_{-\rho}$, there exists a 
slowly 
varying function $L$ such that for all $x >0$, $1-\psi(1/x)=x^{-\rho}L(x)$. Given that the 
distribution function $\psi(1/\cdot)$ 
 is in the domain of attraction of $\Phi_\rho$, the Poisson approximation 
implies that there exists a sequence $(a_n)$ of positive constants such that, for all $x > 0$, 
\begin{equation}\label{eq:zB2}
\tag{C.1}
\lim_{n \to \infty} n \left[1-\psi \{1/(a_nx)\}\right]=  \lim_{n \to \infty} n (a_n x)^{-\rho} 
L(a_n x) =  
x^{-\rho}.
\end{equation}
Furthermore, by Equation~(A6) in the proof of Theorem~2 (a) in \cite{Larsson/Neslehova:2011}, one 
has, for any $j=1, \dots, 
\bar 
\alpha -2$, 
\begin{equation}\label{eq:zB3}
\lim_{x \to \infty} \frac{(-1)^j x^{-j} \psi^{(j)}(1/x)}{\kappa_j x^{-\rho} L(x)} =1, 
\tag{C.2}
\end{equation}
where $\kappa_j = \rho \Gamma(j-\rho)/\Gamma(1-\rho)$. Now for all $i=1, \dots, d$, \Cref{eq:2} 
yields, for any $x > 0$, 
\begin{align*}
      n\Pr\left( \frac{1}{X_i} > a_nx\right) & = n  \left\{1-\bar{H}_i 
\left(\frac{1}{a_nx}\right)\right\}=n 
(a_nx)^{-\rho}L(a_nx) \left\{1-\sum_{j=1}^{\alpha_i-1}\frac{(-1)^j 
(a_nx)^{-j}\psi^{(j)}(1/a_nx)}{j! (a_nx)^{-\rho}L(a_nx)}\right\}.
\end{align*}
Given that $a_n \to \infty$ as $n\to \infty$, the last expression converges by Equations 
\eqref{eq:zB2} and \eqref{eq:zB3} as 
$n\to \infty$ to 
\begin{align*}
x^{-\rho} \left\{ 1- \sum_{j=1}^{\alpha_i -1} \frac{\kappa_j}{j!} \right\} = x^{-\rho} \left\{ 1- 
\rho \sum_{j=1}^{\alpha_i -1} 
\frac{\Gamma(j-\rho)}{\Gamma(j+1)\Gamma(1-\rho)}\right\} = x^{-\rho} \frac{c(\alpha_i, 
-\rho)}{\Gamma(1-\rho)}.
\end{align*}
The Poisson approximation thus implies that, as $n\to \infty$, for all $i=1, \dots, d$ and $x > 0$, 
\begin{equation}\label{eq:zB4}
      \bar{H}_i^n\left( \frac{1}{a_nx} \right)\to \exp 
\left\{-x^{-\rho}\frac{c(\alpha_i, -\rho)}{\Gamma(1-\rho)} 
\right\}.\tag{C.3}
\end{equation} 
For any $\bs{x} \in (0, \infty)^d$, let $1/(a_n \bs{x}) = \{1/(a_n x_1), \dots, 1/(a_n x_d)\}$ and 
denote by $\bar{{x}}_H$ 
the 
harmonic mean of $\bs{x}$, viz. $\bar{x}_H = d/(1/x_1 + \cdots + 1/x_d)$. From  \Cref{eq:1} one 
then has
\begin{align*}
n  \left\{1-\bar{H}\left(\frac{1}{a_n\bs{x}}\right)\right\} = n \left( \frac{a_n \bar{x}_H}{d}  
\right)^{-\rho} L\left( 
\frac{a_n \bar{x}_H}{d}  \right)
\left\{ 1- \sum_{\substack{(j_1, \dots, j_d ) \in \mathbb{I}_{\bs{\alpha}} \\ (j_1, \dots, j_d) 
\neq \bs{0}_d}} 
\frac{(-1)^{j_1+ 
\cdots + j_d} \left( \frac{a_n \bar{x}_H}{d}  \right)^{-j_1 - \dots - j_d}\psi^{(j_1+ \cdots + 
j_d)} \left( \frac{d}{a_n 
\bar{x}_H} \right) }{j_1! \dotsm j_d! \left( \frac{a_n \bar{x}_H}{d}  \right)^{-\rho} L\left( 
\frac{a_n 
\bar{x}_H}{d}  \right)} \prod_{i=1}^d \left( \frac{\bar{x}_H}{d x_i} \right)^{j_i} \right\}
\end{align*}
By \Cref{eq:zB3}, the right most expression in the curly brackets converges, as $n\to \infty$, to
\begin{multline*}
1- \rho\sum_{\substack{(j_1, \dots, j_d ) \in \mathbb{I}_{\bs{\alpha}} \\ (j_1, \dots, j_d) \neq 
\bs{0}_d}} 
\frac{\Gamma(j_1 + \cdots + j_d - \rho)}{\Gamma(1-\rho) j_1! \dotsm j_d!} \prod_{i=1}^d \left( 
\frac{\bar{x}_H}{d x_i} 
\right)^{j_i}  \\ 
=   1- \rho\sum_{\substack{(j_1, \dots, j_d ) \in \mathbb{I}_{\bs{\alpha}} \\ (j_1, \dots, j_d) 
\neq \bs{0}_d}} 
\frac{\Gamma(j_1 + \cdots + j_d - \rho)}{\Gamma(1-\rho)} \prod_{i=1}^d 
\frac{1}{\Gamma(j_i+1)}\left( \frac{1/x_i}{1/x_1 + \dotsm 
+ 1/x_d} \right)^{j_i}.
\end{multline*}
Furthermore, \Cref{eq:zB2} implies that, as $n \to \infty$, 
\begin{align*}
n \left( \frac{a_n \bar{x}_H}{d}  \right)^{-\rho} L\left( \frac{a_n \bar{x}_H}{d}  \right) \to 
\left( \frac{1}{x_1} + 
\cdots + \frac{1}{x_d}\right)^{\rho}.
\end{align*}
Consequently, as $n\to \infty$, $ n \{1-\bar{H}(1/a_n\bs{x})\}  \to -\log H_0(\bs{x})$, where 
\begin{align*}
-\log H_0(\bs{x}) = \left( \frac{1}{x_1} + \cdots + \frac{1}{x_d}\right)^{\rho}   \left\{ 1- 
\rho\quad\sum_{\mathclap{\substack{(j_1, \dots, j_d ) \in \mathbb{I}_{\bs{\alpha}} \\ (j_1, \dots, 
j_d) \neq \bs{0}_d}}}\quad 
\frac{\Gamma(j_1 + \cdots + j_d - \rho)}{\Gamma(1-\rho)} \prod_{i=1}^d 
\frac{1}{\Gamma(j_i+1)}\left( 
\frac{1/x_i}{\sum_{j=1}^d \frac{1}{x_j}
} \right)^{j_i} \right\}.
\end{align*}
By \Cref{eq:B2}, $1/\bs{X} \in \mathcal{M}(H_0)$.
From \Cref{eq:zB4}, the univariate margins of $H_0$ are scaled Fr\'echet, and Sklar's theorem 
implies that the 
unique copula of $H_0$ is of 
the form \eqref{eq:evc} with stable tail dependence function as  in \Cref{prop:3}.
\qed

\bigskip
\noindent
{ \noindent \it Proof of \Cref{prop:4}.} In view of \Cref{cor:2} and \Cref{thm:1} in 
\cite{Larsson/Neslehova:2011}, it only remains to 
compute the expression for $\ell^{\mathrm{pD}}$ given in part (a). Suppose that $\psi \in 
\mathcal{R}_{-\rho}$ for some $\rho 
> 
0$. This means that there exists a slowly varying function such that for all $x > 0$, $\psi(x) = 
x^{-\rho} L(x)$. Because 
$\psi$ 
is itself a survival function, $\psi \in \mathcal{M}(\Phi_\rho)$ and by the univariate Poisson 
approximation, there exists a 
sequence $(a_n)$ of strictly positive constants such that, for all $x > 0$, 
\begin{align}\label{eq:C4}
\lim_{n\to \infty} n \psi(a_n x) = x^{-\rho}.
\tag{C.4}
\end{align}
Furthermore, by Equation~(A1) in the proof of Theorem~1 (a) in \cite{Larsson/Neslehova:2011}, one 
has, for any $j=1, \dots, 
\bar 
\alpha -1$, 
\begin{align}\label{eq:C5}
\lim_{x\to \infty} \frac{(-1)^j x^j \psi^{(j)}(x)}{\psi(x)} =c(j, \rho).
\tag{C.5}
\end{align}
Now let $\bs{X}$ be the Dirichlet random vector with parameters $\bs{\alpha}$ and radial part $R$ 
whose Williamson $\bar 
\alpha$-transform is $\psi$. Denote the distribution function of $\bs{X}$ by $H$ and its univariate 
margins by 
$F_i$, $i=1, \dots, d$. 
Then 
for all $i=1, \dots, d$, Equations \eqref{eq:C4} and \eqref{eq:C5} imply that
\begin{align}\label{eq:C6}
\lim_{n \to \infty} n \bar F_i(a_n x) = \lim_{n \to \infty} n \sum_{j = 0}^{\alpha_i -1} 
\frac{(-1)^j (a_nx)^j\psi^{(j)}(a_nx) 
}{j!}  = x^{-\rho}\!\sum_{j=0}^{\alpha_i - 1} \frac{\Gamma(j + \rho)}{\Gamma(\rho) \Gamma(j+1)} = 
\frac{ x^{-\rho}c(\alpha_i, 
\rho)}{ \Gamma(\rho+1)}
\tag{C.6}
\end{align}
and hence, by the Poisson approximation, $F_i^n(x) \to \exp\{-x^{-\rho} c(\alpha_i, \rho)/ 
\Gamma(\rho+1)\}$ as $n\to 
\infty$. 

Next, for arbitrary $k=1, \dots, d$ and $1 \le i_1 < \dots < i_k \le d$, let 
$\mathbb{I}_{(\alpha_{i_1}, \dots, \alpha_{i_k})} 
= 
\{   0, \dots, \alpha_{i_1} -1\} \times \dots \times \{0, \dots, \alpha_{i_k} -1\}$. For any 
$\bs{x} \in (0, \infty)^d$, 
Equations 
\eqref{eq:1}, \eqref{eq:C4} and \eqref{eq:C5} imply that
\begin{align*}
\lim_{n\to \infty}n \Pr(X_{i_1} > x_{i_1}, \dots, X_{i_k} > x_{i_k}) &
 = \lim_{n\to \infty} n \sum_{(j_1, \dots, j_k) \in \mathbb{I}_{(\alpha_{i_1}, \dots, 
\alpha_{i_k})}} (-1)^{j_1 + \dotsm + j_k} 
\frac{\psi^{(j_1 + \cdots + j_k)}\{a_n(x_{i_1} + \cdots + x_{i_k})\}}{j_1!\dotsm j_k!} 
\prod_{m=1}^k (a_n x_{i_m})^{j_m} \\
& = (x_{i_1} + \cdots + x_{i_k})^{-\rho} \sum_{(j_1, \dots, j_k) \in \mathbb{I}_{(\alpha_{i_1}, 
\dots, \alpha_{i_k})}} 
\frac{\Gamma(j_1 + \cdots + j_k + \rho)}{\Gamma(\rho) j_1! \dotsm j_k!} \prod_{m=1}^k \left( 
\frac{x_{i_m}}{x_{i_1} + \cdots + 
x_{i_k}}\right)^{j_m}.
\end{align*}
Therefore, for any $\bs{x} \in (0, \infty)^d$, 
\begin{align*}
\lim_{n\to \infty} n \left\{ \sum_{k=1}^d \sum_{1 \le i_1 < \dots <  i_k\le d} 
(-1)^{k+1}\Pr(X_{i_1} > a_{n} x_{i_1}, \dots, 
X_{i_k} > a_{n} x_{i_k})\right\} = - \log H_0(\bs{x}), 
\end{align*}
where
\begin{align*}
- \log H_0(\bs{x}) = \sum_{k=1}^d \sum_{1 \le i_1 < \dots <  i_k\le d} (-1)^{k+1} 
(x_{i_1} + \cdots + x_{i_k})^{-\rho} \sum_{(j_1, \dots, j_k) \in \mathbb{I}_{(\alpha_{i_1}, \dots, 
\alpha_{i_k})}} 
\frac{\Gamma(j_1 +   \cdots + j_k + \rho)}{\Gamma(\rho) j_1! \dotsm j_k!} \prod_{m=1}^k \left( 
\frac{x_{i_m}}{x_{i_1} + \cdots + 
x_{i_k}}\right)^{j_m}.
\end{align*}
By \Cref{eq:B2}, $\bs{X} \in \mathcal{M}(H_0)$. As argued above, the univariate margins of $H_0$ 
are given, for all 
$i=1, \dots, d$ and $x>0$, by $\exp\{-x^{-\rho} c(\alpha_i, \rho)/ \Gamma(\rho+1)\}$. Sklar's 
theorem thus implies that the 
unique copula of $H_0$ is of the form \eqref{eq:evc} with stable tail dependence function indeed as 
given by the expression in 
part (a).
\qed
\section{Proofs from Section \ref{sec:5}}  \label{app:D}

{ \noindent \it Proof of \Cref{prop:5}.} First, we show that for any $\rho >-\min(\alpha_1, 
\dots, \alpha_d)$, 
$\rho \neq 0$, 
\begin{align}\label{eq:D1}
\MoveEqLeft \e\left[\max_{1 \le i \le d} \left\{\frac{x_i D_i^\rho}{c(\alpha_i, 
\rho)}\right\}\right]
= \frac{\Gamma(\bar \alpha)}{|\rho|^{d-1} \prod_{i=1}^d \Gamma(\alpha_i )}\int_{\mathbb{S}_d} 
\max(x_i t_i) \left[ 
\sum_{i=1}^d \{c(\alpha_i, \rho) t_i\}^{1/\rho}\right]^{-\rho-\bar \alpha} \prod_{i=1}^d 
\{c(\alpha_i, 
\rho)\}^{\alpha_i/\rho} (t_i)^{\alpha_i/\rho-1} \d \bs{t}.
\tag{D.1}
\end{align}
Indeed, using the fact that $(D_1, \dots, D_d) \eqdis  \bs{Z} / \|\bs{Z}\|$, where $Z_i \sim 
\mathsf{Ga}(\alpha_i, 1)$, 
$i=1, \dots, d$ are independent, 
\begin{align*}
\e\left[\max_{1 \le i \le d} \left\{\frac{x_i D_i^\rho}{c(\alpha_i, \rho)}\right\}\right] = 
\int_{\mathbb{R}_+^d} \max_{1 
\le   i \le d} \left\{\frac{x_i z_i^\rho}{c(\alpha_i, \rho)}\right\} (z_1+ \cdots + z_d)^{-\rho} 
\prod_{i=1}^d 
\frac{e^{-z_i}   z_i^{\alpha_i-1}}{\Gamma(\alpha_i)} \d \bs{z}.
\end{align*}
Make a change of variable $t_i =\{z_i^\rho/c(\alpha_i, \rho)\}/ \sum_{j=1}^d z_j^\rho/c(\alpha_j, 
\rho)$ for 
$i=1, \dots, d-1$ and $w=\sum_{j=1}^d z_j^\rho/c(\alpha_j, \rho)$. For ease of notation, set also 
$t_d = 1-\sum_{i=1}^{d-1} 
t_i$. 
Then, for $i=1, \dots, d$, $z_i = \{c(\alpha_i, \rho)t_i w \}^{1/\rho}$ and the absolute value of 
the Jacobian is
\begin{align*}
|\mathbf{J}| = \frac{1}{|\rho|^d} w^{d/\rho  -1} \prod_{i=1}^d c(\alpha_i, \rho)^{1/\rho} 
t_i^{1/\rho -1}.
\end{align*}
Therefore, 
\begin{multline*}
\textrm{E}\left[\max_{1 \le i \le d} \left\{\frac{x_i D_i^\rho}{c(\alpha_i, \rho)}\right\}\right] = 
\frac{1}{|\rho|^d\prod_{i=1}^d \Gamma(\alpha_i)}\int_{\mathbb{S}_d} \max_{1 \le i \le d} (x_i t_i) 
\left[ \sum_{i=1}^d 
\{c(\alpha_i, \rho) t_i\}^{1/\rho}\right]^{-\rho} \prod_{i=1}^d c(\alpha_i, \rho)^{\alpha_i/\rho} 
t_i^{\alpha_i/\rho 
-1} \\
\times \int_0^\infty w^{\bar\alpha/\rho-1} e^{-w^{1/\rho}\sum_{i=1}^d \{c(\alpha_i, \rho) 
t_i\}^{1/\rho}} \d w \d 
\bs{t}.
\end{multline*}
\Cref{eq:D1} now follows from the fact that
\begin{align*}
\int_0^\infty w^{\bar\alpha/\rho-1} e^{-w^{1/\rho}\sum_{i=1}^d \{c(\alpha_i, \rho) t_i\}^{1/\rho}} 
\d w = |\rho|\, 
\Gamma(\bar \alpha)\left[\sum_{i=1}^d \{c(\alpha_i, \rho) t_i\}^{1/\rho}\right]^{-\bar 
\alpha}\!\!\!.
\end{align*}
The expression for $h^{\mathrm{D}}$ now follows directly from Eqs.~\eqref{eq:stdfspecdens} and 
\eqref{eq:D1}, while the formulas 
for $h^{\mathrm{pD}}$ and $h^{\mathrm{nD}}$ obtain upon setting $\rho=\rho$ and $\rho=-\rho$, 
respectively.
\qed

\section{Proofs from Section \ref{sec:7}} \label{app:E}

{ \noindent \it Proof of \Cref{prop:6}.} 
For $k=1, \dots, d$, the formula for the $k$th order mixed partial derivatives of 
$\ell^{\mathrm{D}}(1/\bs{x})$ can be 
established from \cref{eq:deHpp}. Indeed, if $\bs{V}$ denotes a random vector with 
independent scaled Gamma components $V_i \sim 
\mathsf{sGa}\{1/c(\alpha_i, \rho), 1/\rho, \alpha_i\}$, then the point process representation 
\cref{eq:deHpp} implies that, for all $\bs{x} \in \mathbb{R}_+^d$, 
\begin{align}\label{eq:integratedintens}
\ell^{\mathrm{D}}(1/\bs{x}) = \int_0^\infty \Pr\left(\frac{V_i}{t} > x_i \; \text{for at least 
one 
$i\in \{1, \dots, d\}$} \, \right) \d t = \int_0^\infty \left[1- \prod_{i=1}^d F\left\{x_i 
t; \frac{1}{c(\alpha_i, \rho)}, \frac{1}{\rho}, \alpha_i\right\}\right] \d t.
\end{align}
For any $k=1, \dots, d$, the expression on the right-hand side of 
\cref{eq:integratedintens}  can be differentiated with respect to $x_1, \dots, x_k$ under 
the integral sign. This gives the formulas for $\partial \ell^{\mathrm{D}}(1/\bs{x})/\partial 
x_1 \dots \partial x_k$. When $k=d$, \cref{eq:integratedintens} implies that
\begin{align*}
\frac{\partial^d \ell^{\mathrm{D}}(1/\bs{x})}{\partial x_1 \cdots \partial x_d} & = - 
\int_0^\infty 
t^d \prod_{i=1}^d f\left\{x_i t; \frac{1}{c(\alpha_i, \rho)}, \frac{1}{\rho}, 
\alpha_i\right\}\d t\\
&= \frac{1}{\rho^d} \prod_{j=1}^d\frac{c(\alpha_j, \rho)\{c(\alpha_j, 
\rho)x_j\}^{\alpha_j/\rho-1}}{\Gamma(\alpha_j)}\int_0^\infty t^{\alpha_j/\rho} 
\exp\left[-t^{1/\rho}\sum_{j=1}^d \{c(\alpha_j, \rho)x_j\}^{1/\rho}\right]\d t 
\\&=  \frac{\Gamma(\bar{\alpha}+\rho)}{\rho^{d-1}\left[\sum_{j=1}^d 
\{c(\alpha_j, \rho)x_j\}^{1/\rho}\right]^{\bar{\alpha}+\rho}} \prod_{j=1}^d 
\frac{c(\alpha_j, \rho)\{c(\alpha_j, \rho)x_j\}^{\alpha_j/\rho-1}}{\Gamma(\alpha_j)}, 
\end{align*}
where the last equality follows upon making the change of variable $u=\sum_{j=1}^d \{c(\alpha_j, \rho)x_j\}^{1/\rho} t^{1/\rho}$. 
Alternatively, Theorem~1 in \cite{Coles:1991} implies that that the $d$th order mixed partial derivative of $\ell^\textrm{
D}(1/\bs{x})$ equals $- d \|\bs{x}\|^{-d-1} h^{\mathrm{D}}(\bs{x}/\|\bs{x}\|;\rho, \bs{\alpha})$, which indeed simplifies to $- d 
 h^{\mathrm{D}}(\bs{x};\rho, \bs{\alpha})$ given that $h^{\mathrm{D}}(\bs{x}/\|\bs{x}\|;\rho, \bs{\alpha}) = \| \bs{x}\|^{d+1} 
h^{\mathrm{D}}\left(\bs{x};\rho, \bs{\alpha}\right)$.  

Finally, the formulas for $F(x; a, b, c)$ follow immediately from the fact that the scaled Gamma distribution is also the  
distribution of the random variable $aZ^{1/b}$, where $Z$ is Gamma with shape $c$ and unit scaling.
\qed

\bigskip
{\it Derivation of the gradient score.}
Straightforward calculations show that 
\begin{align*}
\frac{\partial \log d h^{\mathrm{D}}(\bs{x})}{\partial x_i} &=  
-\frac{(\bar{\alpha}+\rho)c(\alpha_i, \rho)^{1/\rho}x_i^{1/\rho-1}}{\rho \sum_{j=1}^d \{c(\alpha_j, \rho)x_j\}^{1/\rho}} + 
\left(\frac{\alpha_i}{\rho}-1\right)\frac{1}{x_i}\\
\frac{\partial^2 \log d h^{\mathrm{D}}(\bs{x})}{\partial x_i \partial x_k} &=  
-\frac{(\bar{\alpha}+\rho)c(\alpha_i, \rho)^{1/\rho}x_i^{1/\rho-1}}{\rho \sum_{j=1}^d \{c(\alpha_j, \rho)x_j\}^{1/\rho}}\left[ 
\left(\frac{1}{\rho}-1\right)\frac{\mathrm{I}_{ik}}{x_i}  - \frac{c(\alpha_k, \rho)^{1/\rho}x_k^{1/\rho-1}}{\rho \sum_{j=1}^d 
\{c(\alpha_j, \rho)x_j\}^{1/\rho}}\right]-\left(\frac{\alpha_i}{\rho}-1\right) 
\frac{\mathrm{I}_{ik}}{x_i^2}.
\end{align*}
\end{appendices}
\end{document}